

\magnification=1200
\settabs 18 \columns


\baselineskip=17 pt

\def\overleftarrow#1{\vbox{\ialign{##\crcr
  \leftarrowfill\crcr\noalign{\kern-1pt\nointerlineskip}
  $\hfil\displaystyle{#1}\hfil$\crcr}}}
\def\overrightarrow#1{\vbox{\ialign{##\crcr
  \leftarrowfill\crcr\noalign{\kern-1pt\nointerlineskip}
  $\hfil\displaystyle{#1}\hfil$\crcr}}}


\def\Crm{\hskip0.5mm \hbox{\rm l\hskip -5.5pt C\/}}

\def\s{\smallskip}

\def\b{\bigskip}
\def\bb{\bigskip\bigskip}
\def\bbb{\bigskip\bigskip\bigskip}

\def\sqr#1#2{{\vcenter{\vbox{\hrule height.#2pt
 \hbox{\vrule width.#2pt height#1pt \kern#1pt
 \vrule width.#2pt} \hrule height.#2pt}}}}

\def\operp{\hbox{${\kern+.25em{\bigcirc}
\kern-.85em\bot\kern+.85em\kern-.25em}$}}
\def\lsim{\;\raise0.3ex\hbox{$<$\kern-0.75em\raise-1.1ex\hbox{$\sim$}}\;}
\def\gsim{\;\raise0.3ex\hbox{$>$\kern-0.75em\raise-1.1ex\hbox{$\sim$}}\;}
\def\no{\noindent}

\def\ce{\centerline}
\def\ve{\vfill\eject}
\def\rdots{\mathinner{\mkern1mu\raise1pt\vbox{\kern7pt\hbox{.}}\mkern2mu
 \raise4pt\hbox{.}\mkern2mu\raise7pt\hbox{.}\mkern1mu}}

\def\e e{$e^+ e^-$ }


\rightline {UCLA/95/TEP/32}
\rightline {August 1995}
\bbb
\ce {\bf GENERALIZATION AND DEFORMATIONS OF QUANTUM GROUPS;}\ce {\bf
QUANTIZATION OF ALL SIMPLE LIE BI-ALGEBRAS}
\b
\ce {C. Fronsdal}
\s
\ce {\it Department of Physics}
\ce {\it University of California, Los Angeles, CA 90095-1547}
\bb

\no {\bf Abstract.}

A large family of ``standard" coboundary Hopf algebras is investigated.
The existence of a universal R-matrix is demonstrated for the case when
the parameters are in general position. Special values of the parameters
are characterized by the appearance of certain ideals; in this case the
universal R-matrix exists on the associated algebraic quotient. In
special cases the quotient is a ``standard" quantum group; all familiar
quantum groups including twisted ones are obtained in this way. In other
special cases one finds new types of coboundary bi-algebras.

A large class of first order deformations of all these standard
bi-algebras is investigated and the associated deformed universal
R-matrices have been calculated. One obtains, in particular, universal
R-matrices associated with all simple, complex Lie algebras
(classification by Belavin and Drinfeld) to first order in the
deformation parameter.

\ve

\line {\bf 1. Introduction.   \hfil}
\s

Quantum groups sprouted in that fertile soil where mathematics overlaps
with physics.  The mathematics of quantum groups is exciting, and the
applications to physical modelling are legion.  It is the more surprising
that some aspects of the structure of quantum groups remain to be
explored; this is especially true of those aspects that bear upon the
problem of classification.  The quantum groups that have so far found
employment in physics are very special (characterized by a single
``deformation" parameter
$q$).  It is true that these applications are susceptible to some
generalization, by the process of twisting;  unfortunately it is easy to
receive the totally erroneous impression that twisting is a gauge
transformation that relates equivalent structures.  The fact that twisted
or multiparameter quantum groups differ qualitatively among themselves
becomes evident when one investigates their rigidity to deformation.
Deformation theory is a means of attacking the classification problem; at
the same time it offers a wider horizon against which to view the whole
subject.  The new quantum groups discovered this way (the deformations of
the twisted ones) are dramatically different; the physical applications
should be of a novel kind.

Let $\cal L$ be a simple Lie algebra over \Crm.  A structure of
coboundary Lie bialgebra on $\cal L$ is determined by a ``classical"
r-matrix; an element $r\in \cal L \otimes \cal L $ that satisfies the
classical Yang-Baxter relation
$$  [r_{12},r_{13} + r_{23}] + [r_{13},r_{23}]=0~, \eqno(1.1)
$$
\no as well as the symmetry condition
$$
 r+r^t=\hat K~, \eqno(1.2)
$$
\no where $\hat K$ is the Killing form of $\cal L$.  The classification
of r-matrices of simple complex Lie algebras was accomplished by Belavin
and Drinfeld [BD].

It is widely believed that there corresponds, to each such r-matrix, via
a process of ``quantization,"  a unique quantum group [D2]. Somewhat more
preciesely, one expects that there exists a Hopf algebra deformation
$\tilde U(\cal L)$ of $U(\cal L)$, and an element
$R\in \tilde U({\cal L})\otimes \tilde U({\cal L})$ such that $\Delta
R=R\Delta^\prime$, where $\Delta$ is the coproduct of $\tilde U(\cal L)$
and
$\Delta^\prime$ is the opposite coproduct,   satisfying the (quantum)
Yang-Baxter relation
$$
 R_{12}R_{13}R_{23}=R_{23}R_{13}R_{12}~; \eqno(1.3)
$$
\no such that $r$ can be recovered by an expansion of $R$ with respect to
a parameter $\hbar$:
$$  R=1+\hbar r+o(\hbar^2)~. \eqno(1.4)
$$

Till now, this program has been realized for r-matrices of a restricted
class that we shall call ``standard".
\vskip.5cm

\no {\bf Definition 1.}  Let $\cal L$ be a simple, complex Lie algebra,
${\cal L}^0$ a Cartan subalgebra and $\Delta^+$ a set of positive roots.
A standard r-matrix for $\cal L$ has the expression
$$
 r=r_0 + \sum_{\alpha\in\Delta^+} E_{-\alpha}\otimes E_\alpha~. \eqno(1.5)
$$
\no Here $r_0\in {\cal L}^0 \otimes {\cal L}^0$ is restricted by (1.2).
\vskip.5cm

The (universal) R-matrix that corresponds to a standard r-matrix is
known. An explicit formula has been given for the simplest choice of
$r_0$ [KR]. Of the more or less explicit formulas for those R-matrices
that correspond to the standard r-matrices in general, there is one that
seems the more fundamental [Ro][LSo]
$$
 R=R^0\bigl(1+\sum_\alpha e_{-\alpha}\otimes e_\alpha + \ldots\bigr)~.
\eqno(1.6)
$$
\no Here $\{H_a,e_\alpha,e_{-\alpha}\}$ are Chevalley-Drinfeld generators
associated with the Cartan subalgebra and simple roots,
$R^0$ involves only the former.  An R-matrix of this form will be called
standard; a precise   definition (in a more general context) will be
given in Section 2, Definition 2.2.  The relationship between (1.5) and
(1.6) will be examined in Section 11.

The R-matrices associated with the twisted quantum groups discovered by
Reshetikhin [R] and others [Sc][Su] are thus all included in the rubrique
``standard".  The principal characteristic of a standard R-matrix is that
it ``commutes with Cartan":
$$
\bigl[H_a\otimes 1+1\otimes H_a,R\bigr]=0.
$$  Non-standard R-matrices are known only in the fundamental
representation [CG][FG].

The original plan of this work was to use deformation theory to obtain
some information about the so far unknown quantum groups that are alleged
to be associated with non-standard r-matrices.  This seems a reasonable
approach because (i) the non-standard r-matrices can be viewed, and
effectively calculated [F], as deformations of standard  r-matrices and
(ii) the largest family of non-standard quantum groups known so far was
found by applying deformation theory to  certain standard R-matrices in
the fundamental representation [FG2].

However, the application of deformation theory to the evaluation of
non-standard   R-matrices turns out to be complicated, perhaps because no
useful cohomological framework could be discovered. If some progress has
been achieved in the present paper, then it is due, in the first place,
to the idea of focusing on the representation (1.6) of the standard
universal R-matrix, and in the second place to the discovery of a
differential complex associated with the Yang-Baxter relation: a detailed
study of (1.6) turned out to be unexpectedly rewarding.

It turns out that the representation (1.6) for an R-matrix that satisfies
the Yang-Baxter equation makes sense in a context that is much wider than
quantum groups.

We introduce an algebra ${\cal{A}}$ (Definition 2.1) with generators
$\{H_a,e_\alpha,e_{-\alpha}\}$ that satisfy certain relations, including
the following (see also Eq. (1.7b) below):
$$
 [H_a,H_b]=0,~[H_a,e_{\pm\alpha}] = \pm H_a(\alpha) e_{\pm\alpha}~,
\eqno(1.7a)
$$
\no with $H_a(\alpha) \in \Crm$.  We define a standard R-matrix on
${\cal{A}}$~--~Definition 2.2.~--~as a formal series of the form
$$
\eqalign{R&=\exp(\varphi^{ab}H_a\otimes H_b)\bigl(1+e_{-\alpha}
\otimes e_\alpha  + \sum^\infty_{k=2} t^{(\alpha^\prime)}_{(\alpha)}
e_{-\alpha_1}\ldots e_{-\alpha_k}\otimes e_{\alpha^\prime_1}\ldots
e_{\alpha^\prime_k}\bigr)~, \cr}
$$
\no with parameters $\varphi^{ab}\in \Crm$, fixed, and try to determine
the coefficients $t^{(\alpha^\prime)}_{(\alpha)}\in  \Crm
 $ so that the Yang-Baxter relation (1.3) is satisfied.  One finds that
this requires additional relations, namely
$$  [e_\alpha,e_{-\beta}] = \delta^\beta_\alpha
\bigl(e^{\phi(\alpha,\cdot)}-e^{-\varphi(\cdot,\alpha)}\bigr)~,
\eqno(1.7b)
$$
\no with $\varphi(\alpha,\cdot)=\varphi^{ab}H_a(\alpha)H_b$,
$\varphi(\cdot,\alpha)=\varphi^{ab}H_aH_b(\alpha)$.  These relations are
therefore included in the definition of the algebra
${\cal{A}}$, Definition (2.1). Generically, with the parameters in
general position, no further relations are required.

The generators $H_a$ of ${\cal{A}}$ generate an Abelian subalgebra that
is denoted ${\cal{A}}^0$ and sometimes called the Cartan subalgebra.  A
key point is to refrain from introducing, \break {\it a priori},   any
(generalized) Serre relations among the  Chevalley-Drinfeld generators
$e_\alpha$, or among the
$e_{-\alpha}$.  The algebras of ultimate interest are obtained
subsequently, by identifying an appropriate ideal $I\subset {\cal{A}}$
that intersects ${\cal{A}}^0$ trivially, and passing to the quotient
${\cal{A}}/I$.  This is the strategy of Chevalley [C], fully exploited in
the theory of Kac-Moody algebras [K][M]; here it is applied to
``generalized quantum groups."

The first result is Theorem 2.  It asserts that the Yang-Baxter relation
for the standard R-matrix on $\cal A$ is equivalent to a simple, linear
recursion relation for the coefficients
$t^{(\alpha^\prime)}_{(\alpha)}$.

The integrability of this recursion relation, Eq.(2.14), is studied in
Sections 3, 4 and 5; it is related to the existence of ``constants," and
eventually to generalized Serre relations.  Generically, there are no
constants and no obstructions, whence the second result that, when the
parameters of ${\cal{A}}$ are in general position, there exists a unique
set of coefficients
$t^{(\alpha^\prime)}_{(\alpha)}$ such that the standard R-matrix
satisfies the Yang-Baxter relation.

Constants exist on certain hyper-surfaces in the space of parameters of
${\cal{A}}$; they represent obstructions to the solution of the recursion
relation (2.14) and thus to the Yang-Baxter relation. Constants are
studied in a slightly more general context in Sections 3 and 4; their
complete classification seems to present a formidable, but not unsolvable
problem.  The relevance of this discussion to Yang-Baxter is demonstrated
in Section 5, and the proof of Theorem 2 can now be completed in Section
6.

The study of the obstructions is taken up again in Section 7.  The third
main result is Theorem 7: the obstructions (that is, the constants)
generate an ideal $I\subset {\cal{A}}$, and a unique standard R-matrix,
satisfying Yang-Baxter, exists on ${\cal{A}}/I$.

In Sections 8-10 we turn to the deformations of the standard R-matrix,
but in a context that is wider than quantum groups.  We calculate a class
of first order deformations of the standard R-matrix on ${\cal{A}}/I$ for
any ideal of obstructions
$I\subset {\cal{A}}$.  This is our fourth result, Theorem 10. The main
difficulty is that the problem is not well posed, for we have been unable
to discover a category that is both natural and convenient in which to
calculate \underbar {all} deformations.  We  limit our study to a class
of deformations.  The good news is Theorem 11: when we specialize   to
the case of simple quantum groups, then we obtain quantizations of all
simple Lie bialgebras (constant r-matrices) of Belavin and Drinfeld.

Finally, in Section 12, we define the coproduct, counit and antipode that
turn all these algebras into Hopf algebras.

\vskip 1in

\line {{\bf 2. Standard Universal R-matrices.} \hfil}
\s

The universal R-matrix of a standard or twisted quantum group has the form
$$
\eqalign{R&=\exp(\varphi^{ab}H_a\otimes H_b)  \cr &\,\,\,\,\times
\bigl(1+t_\alpha(e_{-\alpha}\otimes e_\alpha)+t_{\alpha\beta}
(e_{-\alpha}e_{-\beta}\otimes e_{\alpha}e_{\beta}) +
t^\prime_{\alpha\beta}(e_{-\alpha}e_{-\beta}\otimes e_\beta e_\alpha) +
\ldots \bigr)~. \cr} \eqno(2.1)
$$
\no The $H_a$ are generators of the Cartan subalgebra, the $e_\alpha$ are
generators associated with simple roots,
$\varphi^{ab},~t_\alpha,~
t_{\alpha\beta},~t^\prime_{\alpha\beta},~\ldots$ are in the field; the
unwritten terms are monomials in the
$e_\alpha$ and
$e_{-\alpha}$.

More generally, consider the expression (2.1) in which $H_a,
e_{\pm\alpha}$ generate an associative algebra  with certain relations.
\vskip.5cm

\no {\bf Definition 2.1.}  Let $M,N$ be two countable sets,
$\varphi,\psi$  two maps,
$$
\eqalign{& \varphi :~~M\otimes M \rightarrow \Crm~, \cr & \psi
:~~M\otimes N \rightarrow \Crm~, \cr} \quad
\eqalign{a,b &\rightarrow \varphi^{ab}~, \cr a,\beta & \rightarrow
H_a(\beta)~. \cr} \eqno(2.2)
$$
\vskip.5cm
\no Let ${\cal{A}}$ or ${\cal{A}}(\varphi,\psi)$ be the universal,
associative, unital  algebra over \Crm \enskip with generators
$\{H_a\}\, a\in M,~ \{e_{\pm \alpha}\}\,\alpha \in N$, and relations
$$
\eqalignno{&[H_a,H_b]=0~, \quad  [H_a,e_{\pm\beta}] = \pm
H_a(\beta)e_{\pm\beta}~, & (2.3) \cr
&[e_\alpha,e_{-\beta}]=\delta^\beta_\alpha
\bigl(e^{\varphi(\alpha,\cdot)}-e^{-\varphi(\cdot,\alpha)}\bigr)~, &
(2.4) \cr}
$$
\no with $\varphi(\alpha,\cdot)=\varphi^{ab}H_a(\alpha)H_b,~
\varphi(\cdot,\alpha)=\varphi^{ab}H_aH_b(\alpha)$ and $
e^{\varphi(\alpha, \cdot) + \varphi(\cdot,\alpha)} \neq 1, ~~ \alpha \in N
$.  (The last condition on the parameters is included in order to avoid
having to make some rather trivial exceptions.)

The free subalgebra generated by
$\{e_\alpha\}~\alpha\in N$ (resp.
$\{e_{-\alpha}\}~ \alpha\in N$) will be denoted ${\cal{A}}^+$ (resp.
${\cal{A}}^-$); these subalgebras are $Z^+$-graded, the generators having
grade 1. The subalgebra generated by $\{H_a\}\, a\in M$ is denoted
${\cal{A}}^0$.
\vskip.5cm

\ve

\no {\bf Definition 2.2.}  A standard R-matrix is a formal series of the
form
$$
\eqalign{R &= \exp\bigl(\varphi^{ab}H_a\otimes H_b\bigr)
\bigl(1+e_{-\alpha}\otimes e_\alpha  + \sum^\infty_{k=2}
t^{\alpha^\prime_1\ldots a^\prime_k}_ {\alpha_1\ldots \alpha_k}
e_{-\alpha_1}\ldots e_{-\alpha_k}\otimes e_{\alpha^\prime_1}\ldots
e_{\alpha^\prime_k}\bigr)~. \cr}
\eqno(2.5)
$$
\no In this formula, and in others to follow, summation over repeated
indices is implied.  For fixed
$(\alpha)=\alpha_1,\ldots,\alpha_k$ the sum over $(\alpha^\prime)$ runs
over the permutations of $(\alpha)$. The coefficients
$t^{(\alpha^\prime)}_{(\alpha)}$ are in \Crm.
\vskip.5cm

The special property associated with the qualification ``standard" is
that ``$R$ commutes with Cartan"; indeed
$[R,~H_a\otimes 1+1\otimes H_a]=0,~a\in M$.

We shall determine under what conditions on the parameters
$\varphi^{ab},~H_a(\beta)$ of ${\cal{A}}$, and for what values of the
coefficients $t^{(a^\prime)}_{(\alpha)}$, the R-matrix (2.5) satisfies
the Yang-Baxter relation
$$ Y\hskip-1.0mm B~:=~R_{12}R_{13}R_{23}-R_{23}R_{13}R_{12}=0~. \eqno(2.6)
$$
\no This expression is a formal series in which each term has the form
$\psi_1\otimes \psi_2 \otimes \psi_3\in {\cal{A}}\otimes {\cal{A}}
\otimes {\cal{A}}$.  We assign a double grading as follows.  First extend
the grading of ${\cal{A}}^+$ to the subalgebra of
${\cal{A}}$ that is generated by $\{H_a\} \,a\in M$ and
$\{e_\alpha\}\,\alpha\in N$, by assigning grade zero to $H_a$, and
similarly for
${\cal{A}}^-$. Then $\psi_1$ and $\psi_3$ (but not $\psi_2$) belong to
graded subalgebras of ${\cal{A}}$.  If $\psi_1$ and $\psi_3$ have grades
$\ell$ and $n$, respectively, then define
$$ {\rm grade}~(\psi_1\otimes \psi_2\otimes \psi_3)=(\ell,n)~.
\eqno(2.7)
$$
\no To give a precise meaning to (2.6) we first declare that we mean for
this relation to hold for each grade $(\ell,n)$ separately. This is not
enough, for the number of terms contributing to each  grade is infinite
in general.  The appearance of exponentials in the
$H_a$ can be dealt with in the same way as in the case of simple quantum
groups [V].  If the sets $M,N$ are infinite, then all results are
basis dependent.  Eq.(2.6) means that $Y\hskip-1.0mm B$, projected on any
finite subalgebra of ${\cal A}$, vanishes on each grade; the statement
thus involves only finite sums. The analysis of (2.6) will be organized
by ascending grades.
\vskip.5cm
\no {\bf Remarks.}  (i) It is an immediate consequence of (2.6), in grade
(1,1), that
$$ [e_\alpha,e_{-\beta}]=\delta^\beta_\alpha
\bigl(e^{\varphi(\alpha,\cdot)}-e^{-\varphi(\cdot,\alpha)}\bigr)~.
\eqno(2.8)
$$ This relation was therefore included in the definition of the algebra
${\cal A}$.
\no (ii) No relations of the Serre type have been imposed; in fact no
relations whatever on the subalgebras ${\cal{A}}^+$ and
${\cal{A}}^-$, freely generated respectively by the $e_\alpha$ and the
$e_{-\alpha}$. The contextual meaning of  such relations, including
relations of the Serre type, will be discussed in Sections 3-5 and
especially in Section 7.
\vskip.5cm Before stating the main result, it will be convenient to show
the direct evaluation of $Y\hskip-1.0mm B$ up to grade (2,2).  We expand
$$
 R^0~:=~ \exp(\varphi^{ab} H_a\otimes H_b) = R^i\otimes R_i~, \eqno(2.9)
$$
\no sums over $a,b,i$ implied.  Then
$$  e_{-\alpha}R^i\otimes R_i = R^ie_{-\alpha}\otimes
e^{\varphi(\alpha,\cdot)} R_i~. \eqno(2.10)
$$
\no {\it Grade} (1,1).  The contributions to $R_{12}R_{13}R_{23}$ are of
two kinds:
$$
\eqalign{& R^iR^j e_{-\alpha} \otimes R_iR^k\otimes R_je_\alpha R_k~, \cr
& R^ie_{-\alpha} R^j\otimes R_i e_\alpha R^ke_{-\beta}\otimes
R_jR_ke_\beta~. \cr}
$$
\no Cancellation in $Y\hskip-1.0mm B$ is equivalent to Eq.(2.4).
\vskip 0.5cm
\no {\it Grade} (1,2).  The contributions to $R_{12}R_{13}R_{23}$ are
$$
\eqalign{& R^iR^je_{-\beta}\otimes R_iR^ke_{-\alpha}\otimes R_je_\beta
R_k e_\alpha~, \cr & t^{\alpha'\beta'}_{\alpha\beta}
R^ie_{-\gamma}R^j\otimes R_ie_\gamma R^k e_{-\alpha}e_{-\beta}\otimes
R_jR_ke_{\alpha'} e_{\beta'}~. \cr}
$$
\no Cancellation in $Y\hskip-1.0mm B$ requires that
$$
\eqalign{t^{\alpha\beta}_{\alpha\beta}&=
(1-e^{-\varphi(\alpha,\beta)-\varphi(\beta,\alpha)})^{-1}~, \cr
t^{\beta\alpha}_{\alpha\beta} &= -e^{-\varphi(\beta,\alpha)}
t^{\alpha\beta}_{\alpha\beta}~, \quad \alpha\not=\beta~, \cr
t^{\alpha\alpha}_{\alpha\alpha}&= (1+e^{-\varphi(\alpha,\alpha)})^{-1}~.
\cr}
\eqno(2.11)
$$
\no These conditions are necessary and sufficient that the standard
R-matrix satisfy (2.6) up to grade (2,2).

The obstructions to the existence of coefficients $t^{(\alpha^\prime)}_
{(\alpha)}$ such that $Y\hskip-1.0mm B=0$ up to grade (2,2) are therefore
as follows:
$$
\eqalign{& 1+e^{-\varphi(\alpha,\alpha)}=0~~
\hbox{for some $\alpha\in N$}~, \cr &
1-e^{-\varphi(\alpha,\beta)-\varphi(\beta,\alpha)}=0~~
\hbox{for some pair $\alpha\not= \beta$}~. \cr} \eqno(2.12)
$$
\no They are typical of obstructions encountered at all grades.

Let
$$ t_{\alpha_1\ldots\alpha_l}=
t^{\alpha^\prime_1\ldots\alpha^\prime_\ell}_{\alpha_1\ldots\alpha_\ell}\,
e_{\alpha^\prime_1}\ldots e_{\alpha^\prime_\ell}~. \eqno(2.13)
$$
\vskip.5cm

\no {\bf Theorem 2.}  The standard R-matrix (2.5), on
${\cal{A}}$, satisfies the Yang-Baxter relation (2.6) if and only if the
coefficients $t^{(\alpha^\prime)}_{(\alpha)}$ satisfy the following
recursion relation
$$ [t_{\alpha_1\ldots\alpha_\ell},e_{-\gamma}]=
 e^{\varphi(\gamma,\cdot)}
\delta^\gamma_{\alpha_1}t_{\alpha_2\ldots\alpha_\ell}-
t_{\alpha_1\ldots\alpha_{\ell-1}}\delta^\gamma_{\alpha_\ell}e^{-\varphi(\cdot,\gamma)}~.
\eqno(2.14)
$$
\vskip.5cm

\no {\bf Proof.}  (First part.)  We shall prove that (2.14) is
necessary~--~the ``only if" part.  Then we shall study the integrability
of (2.14).  Later, in Section 6, we shall complete the proof of Theorem 2.
 Insert (2.5) into (2.6) and use (2.10).  The contribution to
$Y\hskip-1.0mm B$ in grade $(\ell,n)$ is
$$
 R_{12}^0R_{13}^0R_{23}^0e_{-\alpha_1}\ldots e_{-\alpha_\ell} \otimes
P^{\gamma_1\ldots\gamma_n}_{\alpha_1\ldots\alpha_\ell} \otimes
e_{\gamma_1}\ldots e_{\gamma_n}~,
$$
\no in which $P$ is the sum over $m$, from 0 to ${\rm min}(\ell,n)$, of
the following elements of ${\cal A}$,
$$
\eqalign{t_{\alpha_{\ell-m+1}\ldots\alpha_\ell}^{\gamma_1\ldots\gamma_m}
&t_{\alpha_1\ldots\alpha_{\ell-m}}e^{-\varphi(\cdot,\sigma)}
t^{\gamma_{m+1}\ldots\gamma_n} \cr
&-t_{\alpha_1\ldots\alpha_m}^{\gamma_{n-m+1}\ldots\gamma_n}
t^{\gamma_1\ldots\gamma_{n-m}} e^{\varphi(\tau,\cdot)}
t_{\alpha_{m+1}\ldots\alpha_\ell}~, \cr} \eqno(2.15)
$$
\no where $\sigma=\gamma_1+\ldots +\gamma_m$ and
$\tau=\alpha_1+\ldots+\alpha_m$.  The Yang-Baxter relation is satisfied
in grade $(\ell,n)$ if and only if this expression, summed over $m$,
vanishes for every index set $\alpha_1,\ldots,\alpha_\ell$ and
$\gamma_1,\ldots,\gamma_n$.  This is so because ${\cal A}^+$  and ${\cal
A}^-$ are freely generated. We have used the definition (2.13) and
$$ t^{\gamma_1\ldots\gamma_n}~:=~
t_{\gamma^\prime_n\ldots\gamma^\prime_n}^{\gamma_1\ldots\gamma_n}
e_{-\gamma^\prime_1}\ldots e_{-\gamma^\prime_n}~. \eqno(2.16)
$$
\no The lowest grades in which $t_\ell=(t_{\alpha_1\ldots\alpha_\ell})$
appears are $(\ell,0)$ and $(0,\ell)$.  In these cases $m=0$ and (2.15)
vanishes identically.  At grade $(\ell,1)$ one finds (summing
$m=0,1$), the linear recursion relation
$$  [t_\ell,e_{-\gamma}]=e^{\varphi(\gamma,\cdot)}
\delta^\gamma_{\alpha_1}t_{\ell-1}-t_{\ell-1}\delta^\gamma_{\alpha_\ell}
e^{-\varphi(\cdot,\gamma)}~, \eqno(2.17)
$$
\no the full expression for which is Eq. (2.14).  This equation is
therefore necessary.  That it is also sufficient will be proved in
Section 6.
\vskip 1in

\line {\bf 3. Differential Algebras. \hfil}
\s

Let $B$ be the unital \Crm-algebra freely generated by
$\{\xi_i\}~i\in  N$, countable.  Suppose given a map
$$ q: N\times N\rightarrow \Crm~, \quad (i,j)\rightarrow q^{ij} \not=0~.
\eqno(3.1)
$$
\no Introduce the natural grading on $B$, $B=\bigoplus B_n$, and a set of
differential operators
$$
\partial_i:~B_n\rightarrow B_{n-1}~,~~i\in N~, \eqno(3.2)
$$
\no defined by
$$
\partial_i\xi_j=\delta_i^j+q^{ij}\xi_j\partial_i~. \eqno(3.3)
$$

We study the problem of integrating sets of equations of the type:
$$
\partial_iX=Y_i~, ~~X\in B~, ~~Y_i\in B~, ~~ i\in N~.   \eqno(3.4)
$$
\no The collection $\{Y_i\}~i\in N$ can be interpreted as the components
of a $B$-valued one-form $Y$, on the space
$\{c^i\partial_i~, ~~ c^i\in \Crm~,~~i\in N\}$.  A constant in
$B_n$ is an element $X\in B_n$, $\partial_iX=0$.

\vskip.5cm

\no {\bf Proposition 3.} (a) The following statements are equivalent: (i)
Eq. (3.4) is integrable for every one-form $Y$ with components in
$B_{n-1}$.  (ii) There are no constants in $B_n$.   (b) When the
parameters $q^{ij}$ are in general position, then there are no constants
in $B_n,~n\geq 1$.
\vskip.5cm

\no {\bf Proof.}  Let
$$ X=X^{i_1\ldots i_n} \xi_{i_1}\dots \xi_n~ \in B_n,
$$
\no then $\partial_iX=0$ means that, for each index set,
$$
\eqalign{X^{i_1\ldots i_n}&+q^{i_1i_2}X^{i_2i_1i_3\ldots}
+q^{i_1i_2}q^{i_1i_3}X^{i_2i_3i_1i_4\ldots} \cr &+\ldots +
q^{i_1i_2}\ldots q^{i_1i_n} X^{i_2\ldots i_ni_1}=0~. \cr} \eqno(3.5)
$$
\no Now fix the unordered index set $\{i_1,\ldots,i_n\}$.  If the values
are distinct then we have a set of $n!$ equations for $n!$ coefficients;
in general the number of equations is always equal to the number of
unknowns.  Solutions exist if and only if the determinant of the matrix
of coefficients vanishes.  This determinant is an algebraic function of
the $q^{ij}$ so that solutions of (3.5), other than $X=0$, exist only on
an algebraic subvariety of parameter space.

The calculation of all these determinants appears to be a formidable
problem.  For $n=2$ the result is
$$D^{12}=1-q^{12}q^{21}~, ~~D^{11}=1+q^{11}~. \eqno(3.6)
$$
\no For $n=3$,
$$
\eqalign{D^{123}&=(1-\sigma^{12})(1-\sigma^{13})(1-\sigma^{23})
(1-\sigma^{12}\sigma^{13}\sigma^{23})~, \cr
D^{112}&=(1+q^{11})(1-\sigma^{12})(1-q^{11}\sigma^{12})~, \cr
D^{111}&=1+q^{11}+(q^{11})^2~, \quad \sigma^{12}~:=~q^{12}q^{21}~. \cr}
\eqno(3.7)
$$
\no It is natural to pass from $B$ to the quotient by the ideal generated
by the constants.  In $B_2$ the constants are
$$
\eqalignno{&\xi_1\xi_2-q^{21}\xi_2\xi_1~, \quad \hbox{when}\quad
\sigma^{12}=1~, & (3.8)\cr &\xi_1\xi_1~~~~~~~~~~~~~~, \quad
\hbox{when}\quad q^{11}=-1~. & (3.9)\cr}
$$
\no If $q^{ii}=-1,~i\in N$ and $\sigma^{ij}=1,~i\not= j$, then the
quotient is a $q$-Grassmann algebra or quantum antiplane.  The constants
in $B_3$ are
$$
\eqalignno{&\xi_1\xi_1\xi_1~~~~~~~~~~~~~~~~~~~~~~~~~~~~~~~~~~~~~~~~~~~~~,
{}~~1+q^{11}+(q^{11})^2=0~, & (3.10)\cr
&\xi_1\xi_1\xi_2-(q^{21})^2~\xi_2\xi_1\xi_1 ~~~~~~~~~~~~~~~~~~~~~~~~~,
{}~~1+q^{11}=0~, & (3.11)\cr
&q^{12}\xi_1\xi_1\xi_2-(1+\sigma^{12})~\xi_1\xi_2\xi_1+
q^{21}\xi_2\xi_1\xi_1~,~~ q^{11}\sigma^{12}=1~; & (3.12)\cr}
$$
\no if $\sigma^{12}=1$, there are two constants
$$
\eqalignno{&q^{11}\xi_1\xi_1\xi_2-(1+q^{11})\xi_1\xi_2\xi_1+
(q^{21})^2~\xi_2\xi_1\xi_1~, & (3.13) \cr
&[[\xi_1,\xi_2]_{q^{21}},\xi_3]_{q^{31}q^{32}}~, \quad
[a,b]_q~:=~ab-qba~, & (3.14) \cr}
$$
\no and finally if $\sigma^{12}\sigma^{13}\sigma^{23}=1$ there is one,
$$
\biggl({1\over q^{31}}-q^{13}\biggr)
(\xi_1\xi_2\xi_3+q^{31}q^{32}q^{21}\xi_3\xi_2\xi_1)+\hbox{cyclic.}
\eqno(3.15)
$$
\no Annulment of (3.8), (3.12) and (3.13) are $q$-deformed Serre
relations [D1].  The last item, Eq. (3.15), may be something new; it
should be interesting to study the quotient of the algebra $B$ with 3
generators by the ideal generated by this constant.

A constant that involves only one variable, $\xi_1$ say, exists if and
only if $q^{11}$ is a root of unity,
$$
\xi_1^n ~~{\rm constant ~~ iff} ~~(q^{11})^n = 1, ~~n = 2,3,... \,.
$$

It is easy to determine all constants of the $q$-Serre type; that is, all
those that involve two geneators and one linearly,
$$ C := \sum^k_{m=0} Q^k_m~(\xi_1)^m\xi_2 (\xi_1)^{k-m}~= 0.
\eqno(3.16)
$$
\no With $q=q^{11}$,
$$
\partial_1(\xi_1)^m=(m)_q~(\xi_1)^{m-1}~, \quad (m)_q~:=~1+q+\ldots
+q^{m-1}~. \eqno(3.17)
$$
\no Setting $\partial_1C=0$ gives, for $q^n\not= 1,~n\in Z\hskip-2mm Z$,
$$ Q^k_m=(-q^{12})^m~q^{m(m-1)/2}
\left(\matrix{k\cr m\cr}\right)_q~, \eqno(3.18)
$$
\no while $\partial_2C=0$ is equivalent to
$$
\prod^{k-1}_{m=0}(1-q^m\sigma^{12})=0~.   \eqno(3.19)
$$
\no When $k=2$, compare $D^{112}$ in Eq. (3.7).  If $k$ is the smallest
integer such that a relation like (3.16) holds, then
$$  1-q^{k-1}\sigma^{12}=0~. \eqno(3.20)
$$

Here are some partial results for $B_4$ and $B_5$. $D_{1234}$ is the
product of 12 factors of the form $1 - \sigma_{ij}$, 4 factors of the
form $1-\sigma_{ij}\sigma_{kl}\sigma_{mn}$, 2 identical factors of the
form $1 - \sigma_{12}\ldots\sigma_{34}$; each group accounts for 24
orders in the $q$'s. $D_{12345}$ is the product of 60 factors of the
first type, 20 factors of the second type, 10 factors of  the third type
and 6 identical factors of the form $1 - $(product of all ten
$\sigma_{ij}, i\neq j)$; each group acconts for 120 orders in the $q$'s.

\ve

\line {\bf 4. Differential Complexes. \hfil}
\s

In the generic case, when there are no constants in $B_n$, the equation
$\partial_iX=Y_i,~Y_i\in B_{n-1},~i\in N$, is always solvable, for any
one-form $Y$.  All one-forms are exact, to be closed has no meaning and
the  differential complex is highly trivial.

The existence of a constant $C\in B_n$ implies that there are one-forms
valued in $B_{n-1}$ that are not exact.  To each 1-dimensional space of
constants in $B_n$ there is a one-dimensional space of non-exact
one-forms, valued in $B_{n-1}$, defined modulo exact one-forms and
obtainable as a limit of $\partial_iX$ as $X\rightarrow C$, after
factoring out a constant.  Thus,
$$ X=\xi_1\xi_2-q^{21}\xi_2\xi_1 \eqno(4.1)
$$
\no becomes a constant as $\sigma^{12}\rightarrow 1$, and a
representative of the associated class of non-exact one-forms is  given by
$$ Y_i=\lim (1-\sigma^{12})^{-1}~\partial_iX =
\cases{\xi_2, & $i=1$, \cr 0~, & $i\not= 1$. \cr} \eqno(4.2)
$$

\vskip 0.5cm
\no {\bf Definition 4.1.} An elementary constant is a linear combination
of re-orderings (permutations of the order of the factors) of a fixed
monomial.
\vskip 0.5cm

A constant $C\in B_n$ also implies a concept of closed one-forms.
\vskip.5cm

\no {\bf Proposition 4.}  If $C\in B_n$,
$$ C=C^{i_1\ldots i_n}\xi_{i_1}\ldots \xi_{i_n}~, \eqno(4.3)
$$
\no is a constant, then the differential operator
$$
\Phi(C)~:=~ C^{i_1\ldots i_n} \partial_{i_1}\ldots\partial_{i_n}
\eqno(4.4)
$$
\no vanishes on $B$.
\vskip.5cm

\no {\bf Proof.}   A constant in $B_n$ is a sum of elementary constants;
it is enough to prove the proposition for the case that $C$ is an
elementary constant.  This implies that there are non-zero $f_i\in \Crm ,
\, i \in N,$ such that the following operator identity
$$
\partial_iC-f_iC\partial_i=0,\,\, i \in N, \eqno(4.5)
$$
\no holds on $B$.  Let $B^*$ be the unital, associative algebra freely
generated by $\{\partial_i\}~i\in N$, and let
$\Phi:~B\rightarrow B^*$ be the unique isomorphism of algebras such that
$\Phi(\xi_i)=\partial_i$. Let $BB^*(q)$ be the unital, associative
algebra generated by
$\{\xi_i,\partial_i\}~i\in N$, with relations (3.3); then $\Phi$ extends
to a unique isomorphism
$$
\Phi^\prime:~BB^*(q)\rightarrow BB^*(\hat q)~, \quad
\hat q^{ij}=1/q^{ji}~.
$$
\no In particular, $\Phi'(\xi_i) = \partial_i$ and
$\Phi^\prime(\partial_i)=-(q^{ii})^{-1}\xi_i,~i\in N$. Now Eq.(4.5) means
that
$\partial_i\circ C=Cf_i\circ
\partial_i$, where $a\circ b$ denotes the product in $BB^*(q)$.  Applying
$\Phi^\prime$ one gets
$$
\Phi(C) \circ \xi_i=(f_i)^{-1}\xi_i\circ \Phi(C)~,
$$
\no implying that $\Phi(C) X=0,~X\in B$.
\vskip.5cm

\no {\bf Definition 4.2.}  Let $C$ be an elementary constant in
$B_n,~n\geq 2$.  A $B_1^*$ one-form $Y$, valued in $B$, will be said to
be $C$-closed if
$$ d_CY~:=~ C^{i_1\ldots i_n}\partial_{i_1}\ldots\partial_{i_{n-1}}
Y_{i_n}=0~.
$$
\vskip.5cm

\no {\bf Examples.}  In $B_2$ the constants are of the type
$C=\xi_1\xi_1$ or $C^\prime=\xi_1\xi_2-q^{21}\xi_2\xi_1$. Now $Y$ is
$C$-closed if $dY~:=~\partial_1Y_1=0$ and $C^\prime$-closed if
$d^\prime Y~:=~\partial_1Y_2-q^{21}\partial_2Y_1=0$.  The first case is
characteristic of Grassmann algebras and the other of quantum planes. Let
${\cal{C}}$ be the collection
$$
\bigl\{C_{ij}=\xi_i\xi_j-q^{ji}\xi_j\xi_i~, \quad i,j\in N~, ~~i\not=
j~\bigr\},
$$
\no and suppose all of them constant.  Then we say that a one-form
$Y$ is ${\cal{C}}$-closed if $Y$ is $C_{ij}$-closed for all
$i\not= j$:
$$
\partial_iY_j-q^{ji}\partial_jY_i=0~, \quad i,j\in N~, ~~i\not= j~.
$$
\no In this case the closure of a $B_1^*$ one-form is expressed in terms
of the $B_1^*$ two-form
$$Z=dY~, \quad Z_{ij}=\partial_iY_j-q^{ji}\partial_jY_i~,
$$
\no and this naturally leads to familiar $q$-deformed de Rham complexes,
with trivial cohomology.
\ve

\line {\bf 5. Integrability of Eq. (2.14). \hfil}
\s

It was seen, in Section 2, that a necessary condition for the standard
R-matrix (2.5) to satisfy the Yang-Baxter relation (2.6) is that the
coefficients $t^{(\alpha^\prime)}_{(\alpha)}$ satisfy (2.14); namely
$$
\eqalign{&[t_{\alpha_1\ldots\alpha_\ell},e_{-\gamma}]=
 e^{\varphi(\gamma,\cdot)}\delta^\gamma_{\alpha_1}
t_{\alpha_2\ldots\alpha_\ell}- t_{\alpha_1\ldots\alpha_{\ell-1}}
\delta^\gamma_{\alpha_\ell}e^{-\varphi(\cdot,\gamma)} ~, \cr &
t_{\alpha_1\ldots\alpha_\ell}~:=~
t^{\alpha^\prime_1\ldots\alpha^\prime_\ell}_ {\alpha_1\ldots\alpha_\ell}
e_{\alpha^\prime_1}\ldots e_{\alpha^\prime_\ell}~. \cr}
\eqno(5.1)
$$
\no Define \footnote*{ This is where we need the last condition in
Definition 2.1.} two differential operators, $\vec\partial\hskip
-1mm_{-\gamma}$  and $\overleftarrow \partial\hskip -1mm_{-\gamma}$, on
${\cal{A}}^+$, by
$$ [X,e_{-\gamma}]=e^{\varphi(\gamma,\cdot)}\vec\partial_{-\gamma}
X-X\overleftarrow\partial\hskip -1mm_{-\gamma}e^{-\varphi(\cdot,\gamma)}~,
\eqno(5.2)
$$
\no $X\in {\cal{A}}^+$; note that $\overleftarrow\partial\hskip
-1mm_{-\gamma}$ operates from the right.  Similarly,
$$ [e_\alpha,Y]=Y\overleftarrow\partial\hskip -1mm_\alpha
e^{\varphi(\alpha,\cdot)}-e^{-\varphi(\cdot,\alpha)}
\vec\partial_\alpha Y \eqno(5.3)
$$
\no defines two differential operators on ${\cal{A}}^-$.  These
definitions are equivalent to the rules
$$
\eqalign{\vec\partial_{-\gamma} e_\alpha &=
\delta^\gamma_\alpha+e^{-\varphi(\gamma,\alpha)} e_\alpha
\vec\partial_{-\gamma}~, \cr  e_\alpha \overleftarrow\partial\hskip
-1mm_{-\gamma} &=
\delta^\gamma_\alpha + e^{-\varphi(\alpha,\gamma)}
\overleftarrow\partial\hskip -1mm_{-\gamma} e_\alpha~, \cr
e_{-\alpha}\overleftarrow\partial\hskip -1mm_\gamma &=
\delta^\gamma_\alpha+e^{-\varphi(\gamma,\alpha)}
\overleftarrow\partial\hskip -1mm_\gamma e_{-\alpha}~, \cr
\vec\partial_\gamma e_{-\alpha} &=
\delta^\gamma_\alpha+e^{-\varphi(\alpha,\gamma)}
e_{-\alpha}\vec\partial_\gamma~. \cr} \eqno(5.4)
$$
\no Eq. (5.1) is equivalent to*
$$
\vec\partial_{-\gamma}t_{\alpha_1\ldots\alpha_\ell}=
\delta^\gamma_{\alpha_1}t_{\alpha_2\ldots\alpha_\ell}~, ~~
t_{\alpha_1\ldots\alpha_\ell}\overleftarrow\partial\hskip -1mm_{-\gamma}=
  t_{\alpha_1\ldots \alpha_{\ell-1}}\delta^\gamma_{\alpha_\ell} ~.
\eqno(5.5)
$$
\vskip.5cm

\no {\bf Proposition 5.}  Suppose that the parameters
$\varphi(\alpha,\beta)$ are in general position, so that there are no
constants in
${\cal{A}}^+~({\cal{A}}^-)$ with respect to the differential operators
$\vec\partial_{-\gamma}$ or
$\overleftarrow\partial\hskip -1mm_{-\gamma}$ ($\vec\partial_\gamma$ or
$\overleftarrow\partial\hskip -1mm_\gamma$).  Then either one of the two
equations in (5.5) determines $t_{\alpha_1\ldots\alpha_\ell}$ recursively
and uniquely (the same in each case) from
$t_\alpha=e_\alpha$.
\vskip.5cm

\no {\bf Proof.}  From (5.4) we deduce that
$$ (\vec\partial_{-\gamma}X)
\overleftarrow\partial\hskip -1mm_{-\gamma^\prime}=
\vec\partial_{-\gamma} (X\overleftarrow\partial\hskip
-1mm_{-\gamma^\prime})~. \eqno(5.6)
$$
\no By Proposition 3, the first of (5.5) determines
$t_\ell=t_{\alpha_1\ldots\alpha_\ell}$ uniquely from
$t_\alpha=e_\alpha$.  The other recursion relation also has a unique
solution, $t_\ell^\prime$ say.  We must show that
$t_\ell=t_\ell^\prime,~\ell > 1$.  Since parentheses are superfluous,
$$
\eqalign{\vec\partial_{-\gamma}\,t_\ell\,
\overleftarrow\partial\hskip -1mm_{-\gamma^\prime} &=
\delta^\gamma_{\alpha_1}\,t_{\ell-1}\,
\overleftarrow\partial\hskip -1mm_{-\gamma^\prime}~, \cr
\vec\partial_{-\gamma} \,t^\prime_\ell\,
\overleftarrow\partial\hskip -1mm_{-\gamma^\prime} &=
\vec\partial_{-\gamma} \,t^\prime_{\ell-1}\,
\delta^{\gamma^\prime}_{\alpha_\ell}~. \cr}
$$
\no Suppose $t^\prime_k=t_k$ for $k=1,\ldots,\ell-1$; then the right-hand
sides are both equal to $\delta^\gamma_{\alpha_1}
\,t_{\ell-2}\,\delta^{\gamma^\prime}_{\alpha_\ell}$.  Then the left-hand
sides are also equal and, since there are no constants,
$t_\ell=t^\prime_\ell$.  Since $t_1=t^\prime_1$
($t_\alpha=t^\prime_\alpha=e_\alpha$), the proposition follows by
induction.
\vskip.5cm We also encounter the relation
$$
\eqalign{&[e_\alpha,t^{\gamma_1\ldots\gamma_n}]=
t^{\gamma_1\ldots\gamma_{n-1}}\delta_{\alpha }^{\gamma_n}
e^{\varphi(\alpha,\cdot)}-
e^{-\varphi(\cdot,\alpha)}\delta^{\gamma_1}_\alpha
t^{\gamma_2\ldots\gamma_n}~, \cr & t^{\gamma_1\ldots\gamma_n}~:=~
t^{\gamma_1\ldots\gamma_n}_{\gamma_1^\prime\ldots\gamma^\prime_n}
e_{-\gamma^\prime_1}\ldots e_{-\gamma^\prime_n}~. \cr} \eqno(5.7)
$$
\no It is equivalent to (5.1) and to either of the following:
$$ t^{\gamma_1\ldots\gamma_n}\overleftarrow\partial\hskip -1mm_\alpha=
 t^{\gamma_1\ldots\gamma_{n-1}} \delta^{\gamma_n}_\alpha~, ~~
\vec\partial_\alpha t^{\gamma_1\ldots\gamma_n}=
\delta_\alpha^{\gamma_1}t^{\gamma_2\ldots\gamma_n}~. \eqno(5.8)
$$
\no The proof is similar to that of Proposition 5.
\vskip.5cm
\vskip 1in

\line {\bf 6. Completion of the Proof of Theorem 2. \hfil}
\s

Suppose that the relations (2.14) are satisfied for $\ell\geq 1$. Now fix
$\ell,n,~\alpha_1,\ldots,\alpha_\ell$ and
$\gamma^1,\ldots,\gamma^n$; we must prove that the expression (2.15),
summed over $m$ from $0$ to ${\rm min}~(\ell,n)$, vanishes.

We begin by calculating the sum over $m=0,1$ (step 1); then we postulate
a formula for the partial sum over $m=0,\ldots,k$ (step $k$). We prove
the formula by induction in $k$, and finally show that the sum vanishes
when $k = {\rm min}~(\ell,n)$.

The term $m=0$ in (2.15) is
$$ [t_\ell,t^n]=
t^{\gamma_1\ldots\gamma_n}_{\gamma_1^\prime\ldots\gamma_n^\prime}
\sum^n_{i=1} e_{-\gamma^\prime_1}\ldots e_{-\gamma^\prime_{i-1}} [t_{\ell
},e_{-\gamma^\prime_i}] e_{-\gamma^\prime_{i+1}}\ldots
e_{-\gamma^\prime_n}~. \eqno(6.1)
$$
\no As in the preceding section we often write $t_\ell,t^n$ for
$t_{\alpha_1\ldots\alpha_\ell}$, $t^{\gamma_1\ldots\gamma_n}$. We shall
gradually make the formulas more schematic so as to bring out their
structure.  By (2.14)
$$
\eqalign{=t^{(\gamma)}_{(\gamma^\prime)} \sum_{i=1}^n
e_{-\gamma^\prime_1}&\ldots \bigl( e^{\varphi(\gamma^\prime_i,\cdot)}
\delta^{\gamma^\prime_i}_{\alpha_1} t_{\ell-1}-t_{\ell-1}
\delta^{\gamma_i^\prime}_{\alpha_\ell} e^{-\varphi(\cdot,\gamma^\prime_i)}
\bigr) \ldots e_{-\gamma^\prime_n}~.\cr} \eqno(6.2)
$$
\no The term $m=1$ is
$$
\eqalign{t^{\gamma_1}_{\alpha_\ell}&t_{\alpha_1\ldots\alpha_{\ell-1}}
e^{-\varphi(\cdot,\gamma_1)}t^{\gamma_2\ldots\gamma_n}
-t^{\gamma_n}_{\alpha_1}t^{\gamma_1\ldots\gamma_{n-1}}
e^{\varphi(\alpha_1,\cdot)}t_{\alpha_2\ldots\alpha_\ell} \cr
&=t_{\ell-1}e^{-\varphi(\cdot,\alpha_\ell)}
\vec \partial_{\alpha_l}t^{\gamma_1\ldots\gamma_n}-
t^{\gamma_1\ldots\gamma_n}\overleftarrow\partial\hskip -1mm_{\alpha_1}
e^{\varphi(\alpha_1,\cdot)}t_{\ell-1} \cr
=t_{\ell-1}t^{(\gamma)}_{(\gamma^\prime)}& \sum
e_{-\gamma^\prime_1}\ldots\delta^{\gamma^\prime_i}_{\alpha_\ell}
e^{-\varphi(\cdot,\gamma^\prime_i)}\ldots e_{-\gamma^\prime_n}
-t^{(\gamma)}_{(\gamma^\prime)} \sum
e_{-\gamma^\prime_1}\ldots\delta^{\gamma^\prime_i}_{\alpha_1}
e^{\varphi(\gamma^\prime_i,\cdot)}
\ldots e_{-\gamma^\prime_n} t_{\ell-1}~. \cr}
$$
\no This agrees with (6.2) except for the position of $t_{\ell-1}$, and
the sign.  Thus, adding the contributions $m=0,1$ we obtain
$$
\eqalign{t^{(\gamma)}_{(\gamma^\prime)} \sum_{i<j}
\bigl\{&e_{-\gamma^\prime_1}\ldots\delta^{\gamma_i^\prime}_{\alpha_1}
e^{\varphi(\gamma^\prime_i,\cdot)}e_{-\gamma^\prime_{i+1}}\ldots
[t_{\ell-1},e_{-\gamma^\prime_j}]\ldots e_{-\gamma^\prime_n} \cr
&+e_{-\gamma^\prime_1}\ldots[t_{\ell-1},e_{-\gamma^\prime_i}]
e_{-\gamma^\prime_{i+1}}\ldots\delta^{\gamma'_j}_{\alpha_l}
e^{-\varphi(\cdot,\gamma^\prime_j)}\ldots e_{-\gamma^\prime_n}\bigr\}~.
\cr}
\eqno(6.3)
$$
\no This completes step 1; all terms involving $t_\ell$ have disappeared
and $t_{\ell-1}$ appears only in commutators that allow us to use (2.14)
again.

We claim that after carrying out step $k$, which includes summing over
$m=0,\ldots,k$, one obtains the following expression
$$
\sum^k_{s=0} \ldots(\delta e^\varphi)^{k-s}\ldots[t_{\ell-k},e_{-\gamma}]
\ldots (\delta e^{-\varphi})^s \ldots~,\,\, k < {\rm min} (l,n),
\eqno(6.4)
$$ and zero, $k = \,\,$min$ (l,n)$. Here the dots stand for products of
the $e_{-\gamma^\prime_i}$, interrupted $k-s$ times by a factor of the
type
$\delta^{\gamma^\prime_i}_ {\alpha_1}
e^{\varphi(\gamma_i^\prime,\cdot)}$, once by $[~,~]$ and $s$ times by a
factor like
$\delta_{\alpha_l}^{\gamma_i}  e^{-\varphi(\cdot,\gamma^\prime_i)}$, as
in (6.3).

To verify this claim we carry out the next step.  We first evaluate the
commutators and examine the cancellations that take place between
successive values of $s$:
$$
\eqalign{\ldots&[t_{\ell-k},e_{-\gamma^\prime_i}] \ldots  (\delta
e^{-\varphi}) \ldots \cr &+ \ldots (\delta e^\varphi) \ldots
[t_{\ell-k},e_{-\gamma^\prime_i}] \dots \cr &\quad =
\ldots(\delta^{\gamma^\prime_i}_\alpha e^{\varphi(\alpha,\cdot)}
t_{\ell-k-1}-t_{\ell-k-1}\delta^{\gamma^\prime_i}_{\alpha^\prime}
e^{-\varphi(\cdot,\alpha^\prime)}) \ldots(\delta e^{-\varphi})
\ldots \cr &\quad \quad \quad  +(\delta e^\varphi)\ldots
(\delta^{\gamma^\prime_j}_\alpha
e^{\varphi(\alpha,\cdot)}t_{\ell-k-1}-t_{\ell-k-1}
\delta^{\gamma^\prime_j}_{\alpha^\prime}
e^{-\varphi(\cdot,\alpha^\prime)}) \ldots~. \cr}
$$
\no The first term in the first line combines with the second term in the
second line to
$$
\ldots (\delta e^\varphi)\ldots[t_{\ell-k-1},e_{-\gamma}] \ldots (\delta
e^{-\varphi}) \ldots~.
$$
\no Successive terms in (6.4) all combine in this way, to reproduce the
same expression with $k$ replaced by $k+1$, except for the fact that
there is no term in the sequence that precedes and collaborates with the
first term and no term that succeeds and collaborates with the last
term.  It remains, therefore, to be proved that the summand
$m=k+1$ in (2.15) precisely supplies the two missing terms.

By (5.8),
$$
\eqalign{& \vec \partial_\beta t^{\gamma_1\ldots\gamma_n} =
\delta^{\gamma_1}_\beta t^{\gamma_2\ldots\gamma_n}~, \cr &
\vec\partial_{\beta_m}\ldots \vec\partial_{\beta_1} t^{\gamma_1\ldots
\gamma_n}=\delta^{\gamma_1}_{\beta_1}\ldots
\delta^{\gamma_m}_{\beta_m} \,t^{\gamma_{m+1}\ldots\gamma_n}~, \cr &
t^{\beta_1\ldots\beta_m}_{\alpha_{\ell-m+1}\ldots\alpha_\ell}
\vec\partial_{\beta_m}\ldots\vec\partial_{\beta_1}
t^{\gamma_1\ldots\gamma_n} =
t^{\gamma_1\ldots\gamma_m}_{\alpha_{\ell-m+1}\ldots\alpha_l}
t^{\gamma_{m+1}\ldots\gamma_n}~. \cr}
$$
\no Hence, if $\vec t_{\alpha_1\ldots\alpha_\ell}$ is the differential
operator
$$
\vec t_{\alpha_1\ldots\alpha_\ell}~:=~
t^{\alpha_1^\prime\ldots\alpha^\prime_\ell}_{\alpha_1\ldots\alpha_\ell}
\vec\partial_{\alpha^\prime_l} \ldots \vec\partial_{\alpha^\prime_1}~,
$$
\no then
$$ t^{\gamma_1\ldots\gamma_m}_{\alpha_{\ell-m+1}\ldots\alpha_\ell}
t^{\gamma_{m+1}\ldots\gamma_n}=
\vec t_{\alpha_\ell\ldots\alpha_{\ell-m+1}} t^{\gamma_1\ldots\gamma_n}
\eqno(6.5)
$$
\no and the first of the two  terms in (2.15) is
$$
\eqalign{t_{\alpha_1\ldots\alpha_{\ell-m}} &e^{-\varphi(\cdot,\sigma)}
\vec t_{\alpha_\ell\ldots\alpha_{\ell-m+1}} t^{\gamma_1\ldots\gamma_n} \cr
&=t_{\alpha_1\ldots\alpha_{\ell-m}}t^{(\gamma)}_{(\gamma^\prime)}
\sum^n_{i=1} e_{-\gamma^\prime_1}\ldots
\bigl[e^{-\varphi(\cdot,\sigma)}
\vec t_{\alpha_\ell\ldots\alpha_{\ell-m+1}}, e_{-\gamma^\prime_i}\bigr]
\ldots \cr  &= t_{\alpha_1\ldots\alpha_{\ell-m}}
t^{(\gamma)}_{(\gamma^\prime)}
\sum_i e_{-\gamma^\prime_1} \ldots
\delta^{\gamma^\prime_i}_{\alpha_{\ell-m+1}}
e^{-\varphi(\cdot,\sigma)}\vec
t_{\alpha_\ell\ldots\alpha_{\ell-m+2}}\ldots~. \cr}
$$
\no By iteration of these steps one finally ends up, when $m=k+1$, with
precisely the missing terms; actually one of the missing terms, we leave
it to the reader to carry out the calculation for the other one.  This
done, the proof of Theorem 2 is complete.
\vskip.5cm

{\bf Corollary 5.} With the parameters in general position, there exists
a unique R-matrix on ${\cal A}$ that satisfies the Yang-Baxter relation.

\line {7. \bf Obstructions and Generalized Serre Relations. \hfil}

We have been concerned with the construction of a standard R-matrix,
Definition (2.2), that satisfies the Yang-Baxter relation,  Eq.(2.6), in
terms of the generators of an algebra
${\cal{A}}$, Definition (2.1). The relations of ${\cal{A}}$ involve
parameters; at certain hypersurfaces in parameter space we have
encountered obstructions, characterized by the vanishing of one or more
of the determinants that we have studied in Section 3.  At these points
there appear elements in ${\cal{A}}^+$ that are constants with respect to
differential operators
$\vec\partial_{-\alpha}$ and/or $\overleftarrow\partial\hskip -1mm
_{-\alpha}$, and elements in ${\cal{A}}^-$ that are constants with
respect to
$\vec\partial_\alpha$ and/or $\overleftarrow\partial\hskip -1mm_\alpha$.

We shall show that all these obstructions can be overcome by the
introduction of additional relations in the definition of ${\cal{A}}$ or,
what is the same, by replacing ${\cal{A}}$ by a quotient ${\cal{A}}/I$,
where $I$ is an ideal generated by the constants.  The next three
propositions relate the null-spaces of the four differential operators to
each other.
\vskip.5cm

\no {\bf Proposition 7.1.}  The space of constants with respect to
$\vec\partial_{-\gamma}$ in ${\cal{A}}^+_n$ has the same dimension as the
space of constants with respect to
$\overleftarrow\partial\hskip -1mm_{-\gamma}$.  If there are no constants
in
${\cal{A}}^+_\ell$ for $\ell=1,\ldots,n-1$, then the two spaces coincide.
\vskip.5cm

\no {\bf Proof.}  An easy consequence of Eq. (5.6).
\vskip.5cm

Let $C\in {\cal{A}}^+_n$ be a constant with respect to
$\vec\partial_{-\gamma}$,
$\gamma \in N$, and assume, provisionally, that there are no constants in
${\cal{A}}^+_\ell,~1<\ell<n$. Without essential loss of generality we
take $C$ to be a linear combination
$$ C=C^{\gamma_1\ldots\gamma_n} e_{\gamma_1}\ldots e_{\gamma_n}~,
$$
\no where the summation runs over the permutations of a fixed set
$\{\gamma_1\ldots\gamma_n\}$, hence over a finite set.  Let $d$ be the
operator that takes $X\in {\cal{A}}^+_n$ to the one-form $Y$ valued in
${\cal{A}}^+_{n-1}$ with components $\vec\partial_{-\gamma}X,~\gamma\in
N$. This operator is represented by a direct sum of finite square
matrices, also denoted $d$.  The constant $C$ is a null-vector for $d$.
The transposed matrix   also has a null-vector; it exists by virtue of
the fact that $dX$ is
$C$-closed: (Definition 4.2):
$$
\eqalign{(C^{\gamma_1\ldots\gamma_n}&
\vec\partial\hskip -1mm_{-\gamma_1}\ldots\vec\partial\hskip
-1mm_{-\gamma_n}) e_{\alpha_1}\ldots e_{\alpha_n} \cr &=
(C^{\gamma_1\ldots\gamma_n}\vec\partial\hskip -1mm_{-\gamma_1}
\ldots\vec\partial\hskip -1mm_{-\gamma_{n-1}})
(d^{\gamma_n\beta_1\ldots\beta_{n-1}}_{\alpha_1\ldots\alpha_n}
e_{\beta_1}\ldots e_{\beta_{n-1}})
\cr &= d^{\gamma_n\beta_1\ldots\beta_{n-1}}_{\alpha_1\ldots\alpha_n}
(C^{\gamma_1\ldots\gamma_n}\vec\partial\hskip -1mm_{-\gamma_1}\ldots
\vec\partial\hskip -1mm_{-\gamma_{n-1}} e_{\beta_1}\ldots
e_{\beta_{n-1}}) =0~. \cr}
$$
\no The obstruction to solving Eq. (5.5) is that the right-hand side is
not in the null-space of the transpose of $d$; it is not
$C$-closed. Indeed, since there are no constants in
${\cal{A}}^+_\ell,~\ell<n$,
$$ (C^{\gamma_1\ldots\gamma_n}\vec\partial\hskip -1mm_{-\gamma_1}\ldots
\vec\partial\hskip -1mm_{-\gamma_{n-1}}) \delta^{\alpha_1}_{\gamma_n}
t_{\alpha_2\ldots\alpha_n}=C^{\alpha_n\ldots\alpha_1}\not= 0~.
$$
\no Recall that the R-matrix is expressed in terms of
$e_{-\alpha_1}\ldots e_{-\alpha_n}t_{\alpha_1\ldots\alpha_n}$. The
obstruction to Yang-Baxter is thus
$$ e_{-\alpha_1}\ldots e_{-\alpha_n}C^{\alpha_n\ldots\alpha_1}~=:~
C^\prime \in {\cal{A}}^-_n~.
$$
\vskip.5cm

\no {\bf Proposition 7.2.}  The element $C^\prime \in {\cal{A}}^-_n$ is a
constant.
\vskip.5cm

\no {\bf Proof.}  One verifies directly that
$\vec\partial_{-\gamma} C=0,~\gamma \in N$, is equivalent to
$C^\prime\overleftarrow\partial\hskip -1mm_\gamma =0,~\gamma\in N$.
\vskip.5cm

Thus, if the first obstruction to the Yang-Baxter relation is
encountered at the evaluation of $t_{\alpha_1\ldots\alpha_n}$, then this
obstruction can be avoided by   replacing ${\cal{A}}$ by the quotient
${\cal{A}}/I_n$, where $I_n$ is the ideal generated by the constants in
${\cal{A}}^\pm_n$.  Once this has been done, we study the obstructions at
the next level.  Since the constants at level
$n$ have been removed there are none in ${\cal{A}}^\pm_\ell,~
\ell\leq n$, and substantially the same analysis applies to constants in
${\cal{A}}^\pm_{n+1}$.  To formulate the final result we need:
\vskip.5cm

\no {\bf Proposition 7.3.}  The ideal $I^+ \subset {\cal{A}}^+$ generated
by the constants of $\vec\partial_{-\gamma}$, coincides with that
generated by the constants of
$\overleftarrow\partial\hskip -1mm_{-\gamma}$.  The same statement holds
true in ${\cal{A}}^-$, {\it mutatis mutandis}.
\vskip.5cm

When the parameters are in general position there are no constants, and
Theorem 2 with Proposition 5 assures us that there is a unique standard
R-matrix in ${\cal{A}}\otimes {\cal{A}}$ that satisfies the Yang-Baxter
relation  (2.6).  We are now in a position to allow for the appearance of
constants.

\vskip.5cm

\no {\bf Remark.}  There are no constants in ${\cal{A}}_1^\pm$; the
generators $H_a,e_{\pm\alpha}$ of ${\cal{A}}$ are also generators of
${\cal{A}}^\prime={\cal{A}}/I$.
\vskip.5cm

\no {\bf Theorem 7.}  Let $I\subset {\cal{A}}$ be the ideal generated by
the constants in ${\cal{A}}^+$ and the constants in
${\cal{A}}^-$, and let ${\cal{A}}^\prime$ be the quotient ${\cal{A}}/I$.
Interpret the standard, universal R-matrix (2.5) as an element of
${\cal{A}}^\prime\otimes {\cal{A}}^\prime$.   The Yang-Baxter relation
for the standard R-matrix on $\cal A'$ is equivalent to the recursion
relation
$$
\eqalign{[t_l,e_{-\gamma}] &= (e_{-\gamma} \otimes
e^{\varphi(\gamma,\cdot)})t_{l-1} - t_{l-1}(e_{-\gamma} \otimes
e^{-\varphi(\cdot,\gamma)}),\cr t_l &:=
t^{(\alpha')}_{(\alpha)}e_{-\alpha_1} \ldots e_{-\alpha_l} \otimes
e_{\alpha_1'}\ldots e_{\alpha_l'}, \cr } \eqno(7.1)
$$ and to either one of the following
$$ [e_\gamma,t_l] = t_{l-1}(e^{\varphi(\gamma,\cdot)} \otimes e_\gamma) -
(e^{-\varphi(\cdot,\gamma)} \otimes e_\gamma)t_{l-1},
\eqno(7.2)
$$
$$ (1 \otimes \vec\partial_{-\gamma})t_l= (e_{-\gamma} \otimes 1)t_{l-1},
\quad t_l(1 \otimes \overleftarrow\partial\hskip -1mm _{-\gamma}) =
t_{l-1}( e_{-\gamma} \otimes 1),\eqno(7.3)
$$
$$ t_l(\overleftarrow\partial\hskip-1mm _\gamma \otimes 1) = t_{l-1}(1
\otimes e_\gamma),\quad (\vec\partial_\gamma \otimes 1)t_l =  (1  \otimes
e_\gamma)t_{l-1}.\eqno(7.4)
$$ These relations are integrable (with $t_1 = e_{-\alpha} \otimes
e_{\alpha}$) and yield a unique standard R-matrix on $\cal A'$.
\vskip1in

\line {\bf 8. Deformations. \hfil}
\s

This work was initiated with the aim of calculating the universal
R-matrices associated with simple Lie algebras, as deformations of the
standard universal R-matrix.  We shall establish a direct correspondence
between the classical r-matrices of Belavin and Drinfeld on the one hand,
and the deformations of the standard, universal R-matrix for simple
quantum groups on the other. In preparation for this we have explored the
meaning of the Yang-Baxter relation in a much more general context, and
we shall endeavor to maintain this generality in our approach to
deformations.  But, as for the types of deformations,   we shall   limit
our study in a way that seems natural in the context of quantum groups.

A deformation of the standard R-matrix is a formal series
$$ R_\epsilon=R+\epsilon R_1+\epsilon^2 R_2 + \ldots~. \eqno(8.1)
$$
\no Here $R$ is given by Eq. (2.5), with the coefficients determined by
the Yang-Baxter relation, and we attempt to find
$R_1,R_2,\ldots$ so that $R_\epsilon$ will satisfy the same relation to
each order in
$\epsilon$.  To make this program precise, we must specify the nature of
the leading term; the remainder should then be more or less unique.

Recall that $R$ ``commutes with Cartan."  An element $Q\in
{\cal{A}}\otimes {\cal{A}}$ is said to have weight $w$ if
$$ [H_a\otimes 1+1\otimes H_a,~Q] = w_a Q~, \quad w_a\in \Crm,~a\in M~.
\eqno(8.2)
$$
\no Thus $R$ has weight zero.  The image of $Q$ by the projection
${\cal{A}}\otimes {\cal{A}}\rightarrow {\cal{A}}^\prime\otimes
{\cal{A}}^\prime$ has the same weight.  We shall suppose that $R_1$ is
homogeneous (has weight).

Recall further that $R$ is driven by the linear term; by virtue of the
Yang-Baxter relation, $R$ is completely determined by the term
$e_{-\alpha}\otimes e_\alpha$.  It is natural to study deformations that
are driven by a similar term, with fixed, non-zero weight:
$$ R_1=S(e_{\pm\sigma} \otimes e_{\pm\rho}) + \ldots~, \eqno(8.3)
$$
\no with $\sigma,\rho$ fixed and the factor $S$ is in ${\cal{A}}^0$. The
unwritten terms are of higher order, in a sense that we must make precise.
\vskip.5cm

\no {\bf Proposition 8.}  The algebra ${\cal{A}}^\prime={\cal{A}}/I$ is
Z\hskip-1.5mm Z-graded, with grade $e_{\pm\alpha}=\pm 1$, grade
$H_a=0$.  An alternative grading is obtained by reversing the sign.
\vskip.5cm

{\bf Proof.} This is a consequence of the fact that the generators of $I$
are homogeneous; ${\cal A}'$ inherits the grading of $\cal A$.

\vskip.5cm

The standard R-matrix is a formal series $\sum_k \psi^-_k\otimes
\psi^+_k,~\psi^\pm_k\in {\cal{A}}^\prime$.  We use the grading of
Proposition 8 in the second space, the alternative grading in the first
space; then grade $\psi^\pm_k=k$ and $R$ is a formal sum of  terms with
grade $(k,k),~k=0,1,2,\ldots~.$  This grading is an extension of that
used previously, made necessary by the appearance of $e_\sigma$ in the
first space and $e_{-\rho}$ in the second.

With the inclusion of (8.3)   the grades  descend to
$(-1,-1)$.  Finally, the unwritten terms in (8.3) is a series by
ascending grades.  The fact that the grades are bounded below is
fundamental.  We claim that $R_\epsilon$, a formal series in
$\epsilon$, each term a formal series in ascending grades, if it
satisfies the Yang-Baxter relation, is completely determined by the
choice of the two generators $e_{\pm\sigma}$ and $e_{\pm\rho}$ in (8.3).

We shall see that the standard R-matrix on ${\cal A}'$, with the
parameters of
${\cal{A}}'$ in general position, is rigid with respect to deformations
of the type (8.3).  We begin our investigation by establishing some
conditions on the parameters that are necessary for the existence of a
deformation.  In the  sections that follow we shall study each of the
four possibilities envisaged by (8.3) separately.  We organize the
contributions to
$$ Y\hskip-1.0mm B_\epsilon~:=~ R_{\epsilon 12}R_{\epsilon 13}R_{\epsilon
23}- R_{\epsilon 23}R_{\epsilon 13}R_{\epsilon 12}
$$
\no in the same way as the contributions to $Y\hskip-1.0mm B$.  A term
$\psi_1\otimes \psi_2\otimes\psi_3$ is said to have grade $(\ell,n)$ if
$\psi_3$ has grade $n$ and $\psi_1$ has alternative grade
$\ell$. We limit ourselves, in this part, to terms linear in $\epsilon$.
\vskip 1in

\line {\bf 9. Deformations of Types $e_{-\sigma}\otimes e_\rho$,
$e_\sigma \otimes e_\rho$ and $e_{-\sigma} \otimes e_{-\rho}$. \hfil}
\s

We suppose that the driving term in $R_1$ is
$$ S(e_{-\sigma}\otimes e_\rho)~, \quad S\in {\cal{A}}^0\otimes
{\cal{A}}^0~.
$$
\no We examine the contributions to $  Y\hskip-1.0mm B_\epsilon$ of order
$\epsilon$.

The lowest grades are (1,0) and (0,1), with contributions
$$
\eqalign{&(Se_{-\sigma}\otimes e_\rho)_{12}~
R^0_{13}R^0_{23}-R^0_{23}R^0_{13}(Se_{-\sigma}\otimes e_\rho)_{12}~, \cr
&R^0_{12}R^0_{13}(Se_{-\sigma}\otimes e_\rho)_{23}- (Se_{-\sigma}\otimes
e_\rho)_{23}R^0_{13}R^0_{12}~, \cr}
$$
\no respectively.  These vanish if and only if
$$
e^{\varphi(\sigma,\cdot)-\varphi(\rho,\cdot)}=1=e^{\varphi(\cdot,\sigma)-\varphi(\cdot,\rho)}~.
\eqno(9.1)
$$
\no In grade (1,1) we have three contributions, the simplest one is
$$
\eqalign{A&= F\bigl(e_{-\sigma}\otimes \{e^{-\varphi(\cdot,\rho)}-
e^{\varphi(\sigma,\cdot)}\}\otimes e_\rho\bigr)~, \cr F &=
R_{12}^0S_{13}R_{23}^0~. \cr}
$$
\no The other contributions to the same grade are, in view of (9.1),
$$
\eqalign{B&= R^iR^je_{-\alpha}\otimes R_ie_\alpha S^ke_{-\sigma}\otimes
R_j e^{\varphi(\alpha,\cdot)}S_ke_\rho \cr &- R^iR^je_{-\alpha}\otimes
R_iS^ke_\alpha e_{-\sigma}\otimes R_jS_ke_\rho~, \cr}
$$
\no and
$$
\eqalign{C &= S^iR^je_{-\sigma}\otimes R^kS_je_\rho e_{-\alpha}\otimes
R_jR_ke_\alpha \cr &- S^iR^je^{-\varphi(\cdot,\alpha)}e_{-\sigma}\otimes
R^ke_{-\alpha}S_je_\rho\otimes R_jR_ke_\alpha~. \cr}
$$
\no No cancellations occur unless
$$ S=R^0~; \eqno(9.2)
$$
\no the form of the Cartan factor $S$ is already fixed.  Now
$$
\eqalign{B&= F\bigl(e_{-\sigma}\otimes \{e^{\varphi(\sigma,\cdot)}-
e^{-\varphi(\cdot,\sigma)}\}\otimes e_\rho\bigr)~, \cr C&=
F\bigl(e_{-\sigma}\otimes \{e^{\varphi(\rho,\cdot)}-
e^{-\varphi(\cdot,\rho)}\} \otimes e_\rho\bigr)~. \cr}
$$
\no The sum $A+B+C$ vanishes if and only if
$$ e^{\varphi(\rho,\cdot)+\varphi(\cdot,\sigma)}=1~. \eqno(9.3)
$$

Let us take stock.  Conditions (9.1)-(9.3) are necessary.   From (9.1)
and (9.3) it follows in particular that
$$
e^{\varphi(\rho,\alpha)+\varphi(\alpha,\rho)}=1=e^{\varphi(\sigma,\alpha)+\varphi(\alpha,\sigma)}~,
\quad \alpha\in N~. \eqno(9.4)
$$
\no These are conditions that are familiar from our investigation of
constants, see Eq.(3.8).  The relations that are thus implied are
$$ e_\rho e_\alpha-e^{\varphi(\rho,\alpha)}e_\alpha e_\rho = 0 = e_\sigma
e_\alpha -e^{\varphi(\sigma,\alpha)}e_\alpha e_\sigma~,
\,\, \forall\alpha \in N.
$$
  This constitutes a high degree of commutativity in ${\cal{A}}^\prime$
and takes us far away from our main interest in simple quantum groups.
We therefore end our investigation of the type $e_{-\sigma}\otimes
e_\rho$ at this point.

Suppose next that the driving term in $R_1$ is
$$ S(e_\sigma\otimes e_\rho)~, \quad S\in {\cal{A}}^0 \otimes
{\cal{A}}^0~.
$$
\no The lowest grade (in $  Y\hskip-1.0mm B_\epsilon$) is (-1,0), with
just one contribution,
$$
\eqalign{&F^\prime \bigl(e_\sigma\otimes e_\rho\otimes
\{e^{-\varphi(\sigma,\cdot)-\varphi(\rho,\cdot)}-1\}\bigr)~, \cr
&F^\prime~:=~S_{12}R^0_{13}R^0_{23}~; \cr}
$$
\no so we need
$$ e^{\varphi(\sigma,\cdot)+\varphi(\rho,\cdot)}=1~. \eqno(9.5)
$$

The next lowest grade is (-1,1), with two contributions.  The simpler one
is
$$ A_1=F\bigl(e_\sigma\otimes
\{e^{-\varphi(\cdot,\rho)}-e^{-\varphi(\sigma,\cdot)}\}\otimes
e_\rho\bigr)~,
$$
\no and the other one
$$
\eqalign{B_1&= S^iR^je_\sigma\otimes S_iR^ke_\rho e_{-\alpha}\otimes
R_jR_ke_\alpha \cr &- S^iR^je^{-\varphi(\cdot,\alpha)}e_\sigma\otimes
R^ke_{-\alpha} S_ie_\rho\otimes R_jR_ke_\alpha~. \cr}
$$
\no If the two terms in $B_1$ are to combine to something that can be
cancelled by $A_1$ we need $S=(KR^i)\otimes R_i,~K\in {\cal{A}}^0$, in
which case
$$ B_1=F\bigl(Ke_\sigma\otimes \{e^{\varphi(\rho,\cdot)}-
e^{-\varphi(\cdot,\rho)}\}\otimes e_\rho\bigr)~.
$$
\no  We need to make $A_1+B_1=0$.  In view of (9.5) it can be
accomplished in two ways:
$$ {\rm (i)}~~K=1~, \quad \hbox{or} \quad {\rm (ii)}~~
e^{\varphi(\rho,\cdot)+\varphi(\cdot,\rho)}=1~.
$$
\no The second solution is very restrictive; it implies that $e_\rho$
quommutes with all the $e_\alpha,~\alpha\in N$;  besides, it is excluded
by the last condition of Definition 2.1.

We therefore abandon (ii) and adopt (i); that is,
$$ S=R^0~.
$$

The next lowest grade is (0,1), with three contributions.  Using (9.5)
they can be reduced to
$$
\eqalign{A_2 &=
F\bigl(\{1-e^{-\varphi(\cdot,\sigma)-\varphi(\cdot,\rho)}\}\otimes
e_\sigma\otimes e_\rho\bigr)~, \cr B_2 &= F\bigl(e_\sigma
e_{-\alpha}\otimes
e^{\varphi(\rho,\alpha)-\varphi(\cdot,\alpha)}e_\rho\otimes e_\alpha \cr
&~~~~-e_{-\alpha}e_\sigma\otimes  e^{\varphi(\alpha,\cdot)}e_\rho\otimes
e_\alpha \bigr)   \cr C_2 &= F\bigl(e_{-\alpha}e_\sigma\otimes
e^{\varphi(\alpha,\rho)-\varphi(\cdot,\rho)}e_\alpha\otimes e_\rho \cr
&~~~~-e_\sigma e_{-\alpha}\otimes e^{-\varphi(\sigma,\cdot)}
e_\alpha\otimes e_\rho\bigr)~. \cr}
$$
\no Cancellation requires very strong conditions on the parameters
 and we shall stop at this point. Similar results are obtained for
deformations of type
$e_{-\sigma}\otimes e_{-\rho}$.

\ve

\line {\bf 10.  Deformations of Type $e_\sigma\otimes e_{-\rho}$. \hfil}
\s

We come to the last case envisaged in Section 8.  Eq. (8.3), when the
driving term in $R_1$ has the form
$$ S(e_\sigma\otimes e_{-\rho})~, \quad S\in {\cal{A}}^0\otimes
{\cal{A}}^0~. \eqno(10.1)
$$
\no This term has grade (-1,-1); it is the only term in $R_1$ with this
grade, the lowest.  The factor $S$, and all other terms in
$R_1$, are completely determined by the Yang-Baxter relation
$  Y\hskip-1.0mm B_\epsilon=0$ to first order in $\epsilon$.  Besides
(10.1) there is in $R_1$ one other term with only two roots, of the form
$$ S^\prime(e_{-\rho}\otimes e_\sigma)~, S^\prime\in {\cal{A}}^0\otimes
{\cal{A}}^0~;
\eqno(10.2)
$$
\no it has grade (1,1).
\vskip.5cm

\no {\bf Theorem 10.}  Let $R$ be the standard R-matrix described in
Theorem 7.  Suppose that $R+\epsilon R_1$ is a first order deformation,
satisfying the Yang-Baxter relation to first order in
$\epsilon$.  Suppose finally that the term of lowest grade in $R_1$ has
the form (10.1); then  the parameters satisfy
$$ e^{\varphi(\cdot,\rho)+\varphi(\sigma,\cdot)}=1~; \eqno(10.3)
$$ and $R_1$ is uniquely determined and has the expression
$$ R_1=  R(Ke_\sigma\otimes Ke_{-\rho}) - (Ke_{-\rho}\otimes Ke_\sigma)R
{}~, \eqno(10.4)
$$
\no with $K~:=~e^{\varphi(\cdot,\rho)}$.
\vskip.5cm

\no {\bf Proof.}  An easy calculation in the lowest grades shows that
(10.3) is necessary and that $S=K\otimes K$, up to a  numerical factor
that we fix once and for all.

Let $R_1^i,~i=1,2,$ be the two summands in (10.4).  The term of order
$\epsilon$ in $ Y\hskip-1.0mm B_\epsilon$ is the sum of the following six
quantities:
$$
\eqalign{A^i_{12} &= (R^i_1)_{12}R_{13}R_{23}- R_{23}R_{13}(R^i_1)_{12}~,
\cr A^i_{13} &= R_{12}(R^i_1)_{13}R_{23}-R_{23}(R^i_1)_{13}R_{23}~, \cr
A^i_{23} &= R_{12}R_{13}(R^i_1)_{23}-(R^i)_{23}R_{13}R_{12}~, \quad
i=1,2~.
\cr}
\eqno(10.5)
$$
\vskip.5cm
\no {\bf Step 1.}  We begin with the term that contains the lowest grade,
(1,1):
$$
\eqalign{&A^1_{13}=R^i[-\alpha]R^j[-\beta]Ke_\sigma\otimes R_i[\alpha]
R^k[-\gamma]\otimes R_j[\beta]Ke_{-\rho}R_k[\gamma]-\ldots~,
\cr &[-\alpha]\otimes [\alpha]~:=~e_{-\alpha_1}\ldots e_{-\alpha_\ell}
\,t^{(\alpha^\prime)}_{(\alpha)}\,e_{\alpha^\prime_1}\ldots
e_{\alpha^\prime_\ell}~. \cr}
\eqno(10.6)
$$
\no A sum over indices and numbers of indices ($\ell~\alpha$'s,
$m~\beta$'s and $n~\gamma$'s) is understood, and $-\ldots$ stands for the
reflected term.  Using the fact that $R$ satisfies $  Y\hskip-1.0mm B=0$
we can convert (10.6) to
$$ A^1_{13}=R^i[-\alpha]R^j[-\beta]Ke_\sigma\otimes R_i[\alpha]
R^kK[-\gamma]\otimes R_j[\beta]KR_k\bigl[e_{-\rho},[\gamma]\bigr]~ +
\ldots ,
\eqno(10.7)
$$
\no where $+\ldots$ stands for a similar expression that contains a
factor $\bigl[e_\sigma,[-\alpha]\bigr]$ in the first space.  We have used
(10.3) and continue to use this relation without comment.
\vskip.5cm

\no {\bf Step 2.}  Evaluate the commutators in (10.7) using (5.1) and
(5.7).  The result
$$
\eqalign{A^1_{13}&=R^i[-\alpha]R^j[-\beta]Ke_\sigma\otimes
R_i[\alpha]R^kK[-\gamma]e_{-\rho}\otimes R_j[\beta]KR_k[\gamma]K^{-1}  +
\ldots \cr} \eqno(10.8)
$$
\no is a sum of four similar expressions.  Note that the evaluation of
the commutators involves a shift in the summation indices
$\ell,m,n$. The lowest grades in (10.8) are (-1,0) and (0,-1).
\vskip.5cm

\no {\bf Step 3.}  Now write down the full expression for
$A^1_{12}+A^1_{23}$; it also contains four similar terms.  Two of them
cancel two of the terms in (10.8); to verify this the relation
$Y\hskip-1.0mm B=0$ must be invoked.
\vskip.5cm

\no {\bf Step 4.}  Combine the remaining two terms from (10.8) with  the
remaining two terms from $A^1_{12}+A^1_{23}$ and verify that
$$
\eqalign{&A^1_{12}+A^1_{13}+A^1_{23}  \cr &=
R^i[-\alpha]KR^j\bigl[e_\sigma,[-\beta]\bigr]\otimes R_i[\alpha]
KR^ke_{-\rho}[-\gamma]\otimes
R_jK[\beta]R_ke^{\varphi(\rho,\cdot)}[\gamma]  + \ldots \, ,\cr}
\eqno(10.9)
$$
\no where $+\ldots$ stands for a term that contains a factor
$\bigl[[\beta],e_{-\rho}\bigr]$ in the third space.
\vskip.5cm

\no {\bf Step 5.}  Evaluate the commutators (second shift of summation
indices)
$$
\eqalign{&= R^i[-\alpha]R^jK[-\beta]K^{-1}\otimes
R^i[\alpha]KR^ke_{-\rho}[-\gamma]\otimes R_jK[\beta]e_\sigma
R_ke^{\varphi(\rho,\cdot)}[\gamma] + \ldots \cr &=:~X_1+X_2+Y_1+Y_2~.
\cr}  \eqno(10.10)
$$ The lowest grades are now (1,0) and (0,1).
\vskip.5cm

\no {\bf Step 6.}  Two of the four terms in (10.10) are:
$$
\eqalign{ X_1 &= R_{12}R_{13}\bigl\{(Ke_{-\rho}
\otimes Ke_\sigma)R\bigr\}_{23}~, \cr Y_2 &=
-R_{23}R_{13}\bigl\{(Ke_{-\rho}\otimes Ke_\sigma)R\bigr\}_{12}~. \cr}
\eqno(10.11)
$$
\no Now add $A^2_{12}+A^2_{23}$ to (10.10) to get
$$
\eqalign{A^1_{12}&+A^1_{13}+A^1_{23}+A^2_{12}+A^2_{23} \cr &= \tilde
X_1+X_2+Y_1+\tilde Y_2~, \cr} \eqno(10.12)
$$
\no where $\tilde X_1$ and $\tilde Y_2$ are obtained from $X_1$ and
$Y_2$ by adding $A^2_{23}$ and $A^2_{12}$.
\vskip.5cm

\no {\bf Step 7.}  Use the relation $Y\hskip-1.0mm B=0$ to modify the
expressions for $\tilde X_1$ and $\tilde Y_2$; then notice that the four
terms in (10.12) can be combined to two,
$$ = R^iK[-\alpha]R^j[-\beta] \otimes
R_iK[e_{-\rho},[\alpha]]R^k[-\gamma] \otimes Ke_\sigma
R_j[\beta]R_k[\gamma] + \ldots,
\eqno(10.13)
$$ where the other term has a factor $[[-\gamma],e_\sigma]$ in the second
space.
\vskip.5cm

{\bf Step 8.} Evaluate the commutators (third shift of summation indices),
$$ = R^iK[-\alpha]e_{-\rho}R^j[-\beta] \otimes
R_jK[\alpha]K^{-1}R^k[-\gamma] \otimes K e_\sigma R_j[\beta]R_k[\gamma] +
\ldots.
\eqno(10.14)
$$ This expression has four terms; the lowest grade is (1.1).
\vskip.5cm

{\bf Step 9.} Two of the four terms in (10.14) cancel each other because
$Y\hskip-1.0mm B = 0$.
\vskip.5cm

{\bf Step 10.} The remaining two terms add up to $-A_{13}^2$.
\vskip.5cm

This completes the verification of the claim that (10.4) defines a first
order deformation of $R$. To complete the proof of Theorem 10 we must
show that this expression (10.4) is unique. This was done by complete
mathematical induction. We omit the details but point out that the key to
the induction processs is visible in steps 2,5 and 8, where the summation
indices are shifted. Theorem 10 is proved.

Let ${\cal P}$ be the collection of pairs $(\sigma,\rho) \in N \otimes N$
such that (10.3) holds; each distinct pair defines a first order
deformation $R + \epsilon R_1^{\sigma,\rho}$ of $R$. Because these
deformations are only first order they generate a linear space
$$ R_1 = \sum_{\sigma,\rho \in {\cal P}} C_{\sigma,\rho}
R_1^{\sigma,\rho},
\eqno(10.15)
$$  with coefficients in the field. The dimension of this space of first
order deformations is zero for parameters in general position. It remains
zero, generically, when the parameters are such that the ideal $I$
generated by the constants is non-empty and
$R$ is  defined on ${\cal A}/I$.

\vskip1in

\line {\bf 11. Classical R-matrices for Simple Lie Algebras. \hfil}
\s

We shall now specialize, by stages, until we arrive at simple quantum
groups, where a confrontation can be made with the list of r-matrices
obtained by Belavin and Drinfeld [BD].

Suppose that $ {\rm Card}(N)~:=~\ell<\infty$.  Suppose next that the
ideal $I$ (generated by the constants of
${\cal{A}}$) is generated by a complete set of Serre relations; that is,
for each pair
$\alpha,\beta\in N,~\alpha\not= \beta$, there is a smallest positive
integer $k_{\alpha\beta}$ such that there is a relation in
${\cal{A}}/I$ of the form
$$
\sum^{k_{\alpha\beta}}_{m=0} Q_m^{(\alpha,\beta)}(e_\alpha)^m
e_\beta(e_\alpha)^{k_{\alpha\beta}-m} =0~, \eqno(11.1)
$$
\no with coefficients $Q^{(\alpha\beta)}_m$ in the field.  The left side,
as an element of ${\cal{A}}^+$, is a constant, and the  penultimate
paragraph of Section 3 applies.  In particular, the relation (3.20)
becomes
$$
e^{\varphi(\alpha,\beta)+\varphi(\beta,\alpha)+(k_{\alpha\beta}-1)\varphi(\alpha,\alpha)}
=1~, \eqno(11.2)
$$
\no We specialize further by supposing that the exponent vanish,
$$
\varphi(\alpha,\beta)+\varphi(\beta,\alpha)=(1-k_{\alpha\beta})
\varphi(\alpha,\alpha)~, \quad \alpha\not= \beta~. \eqno(11.3)
$$ The form $(\cdot,\cdot)$ defined by
$$ (\alpha,\beta)=\varphi(\alpha,\beta)+\varphi(\beta,\alpha) \eqno(11.4)
$$
\no will be called the restricted Killing form, and the $\ell$-by-$\ell$
matrix with components
$$ A_{\alpha\beta}~:=~1-k_{\alpha\beta}= {2(\alpha,\beta)\over
(\alpha,\alpha)} \eqno(11.5)
$$
\no will be called the generalized Cartan matrix; note that it is
symmetrizable. Finally, a suitable restriction on Card($M$) brings us
to   quantum Kac-Moody algebras [K][M].

The semi-classical limit of $R$ is defined by replacing
$$
\eqalign{&\varphi(\cdot,\cdot)\rightarrow \hbar\varphi(\cdot,\cdot)~,
\cr & e_\alpha=cE_\alpha~, \quad e_{-\alpha}=c'  E_{-\alpha}~,
\quad cc^\prime = \hbar~, \quad \alpha\in N~, \cr} \eqno(11.6)
$$
\no and developing the exponentials to first order in $\hbar$.  Then Eq.
(2.4) becomes
$$
\eqalign{[E_\alpha,E_{-\beta}]&=
\delta_{\alpha\beta}\bigl(\varphi(\alpha,\cdot)
+\varphi(\cdot,\alpha)\bigr)  \cr &=:~\delta_{\alpha\beta}H_{(\alpha)}
\varphi(\alpha,\alpha)~. \cr}
\eqno(11.7)
$$
  It follows from (11.7) and (2.3) that
$$ [H_{(\alpha)},E_\beta]=A_{\alpha\beta}E_\beta~, \quad
\alpha,\beta\in N~. \eqno(11.8)
$$ The definition (11.5) of the generalized Cartan matrix implies that
$A_{\alpha\alpha}=2,~\alpha\in N$, that $A_{\alpha\beta}\in
\{0,-1,-2,\ldots\},~\alpha\not= \beta$, and that $A_{\alpha\beta}\not= 0$
implies $A_{\beta\alpha}\not= 0$.  Special cases are affine Lie algebras
and simple Lie algebras.  The latter are characterized by two additional
properties of $(A_{\alpha\beta})$: indecomposability and
${\rm det}(A)>0$.  We may now assume that both hold, and that
$\{H_{(\alpha)},\, \alpha \in N\}$ generates ${\cal A}^0$.

The (classical) r-matrix $r$, associated with the standard R-matrix (2.5)
is defined by
$$  R=1+\hbar r+o(\hbar^2)~. \eqno(11.9)
$$
\no Two terms in r are obvious: $r=\varphi  + \sum E_{-\alpha}
\otimes E_\alpha+{\rm ?}$, with the sum extending over simple roots.
Evaluating the remaining terms is more difficult, but the result is
known.  With a particular normalization of the non-simple roots,
$$
 r=\varphi+\sum_{\alpha\in\Delta^+} E_{-\alpha}\otimes E_\alpha~,
\eqno(11.10)
$$
\no where $\Delta^+$ is the set of positive roots.  This may be called
the standard r-matrix.  It satisfies the classical Yang-Baxter relation
$$
 [r_{12},r_{13} + r_{23}]+[r_{13},r_{23}]=0 \eqno(11.11)
$$
\no and
$$
 r+r^t=\hat K~ ,\eqno(11.12)
$$ the Killing form of ${\cal L}$. In the list of (constant) r-matrices
obtained by Belavin and  Drinfeld [BD], (11.10) is the simplest.  The
quantum groups to which these r-matrices are associated are the twisted
quantum groups of Reshetikhin and others [R][Sc][Su].

To any first order deformation of $R$, there corresponds a first order
deformation of $r$,
$$R_{\epsilon} = 1 + \hbar r_\epsilon + o(\hbar^2), \quad
r_\epsilon=r+\epsilon r_1+o(\epsilon^2)~. \eqno(11.13)
$$
\no Eqs. (10.4) and (10.5) give us
$$ r_1=\sum_{\sigma,\rho\in {\cal{P}}} C_{\sigma,\rho} (E_\sigma \wedge
E_{-\rho})~, \eqno(11.14)
$$
\no where ${\cal{P}}$ is the set of pairs with the property
$$
\varphi(\rho,\cdot)+\varphi(\cdot,\sigma)~:=~0~. \eqno(11.15)
$$

The original work of Belavin and Drinfeld culminates in a list of
constant r-matrices that is complete up to equivalence.  Their results
have recently been re-derived in terms of deformation theory and the
associated cohomology.
\vskip.5cm

\no {\bf Proposition 11.} [F] Let $r$ be the standard r-matrix (11.10)
for a simple Lie algebra ${\cal{L}}$.  The space of essential, first
order deformations of $r$, satisfying (11.11) and (11.12), is
$$ H^2({\cal L}^*,\Crm)=\bigl\{r_1=\sum_{\sigma,\rho\in {\cal{P}}}
C_{\sigma,\rho}E_\sigma \wedge E_{-\rho}+ \sum \tilde C^{ab} H_a\otimes
H_b\bigr\}~. \eqno(11.16)
$$
\no The exact deformations are of finite order and coincide with the
r-matrices of [BD].

\vskip.5cm The second, Cartan term is not ``essential" in the present
context; it represents the freedom to vary the parameters.  We conclude
that
\vskip.5cm

\no {\bf Theorem 11.}  The first order deformations of the standard
R-matrix calculated in Section 10, upon specialization to a simple
quantum group, are in one-to-one correspondence, via (11.9), with the
essential first order deformations of the associated  r-matrix
(11.10)-modulo variations of the parameters.

\vskip 1in

\line {\bf 12. Hopf Structure. \hfil}
\s

It is of some interest to verify that the standard R-matrix, satisfying
the Yang-Baxter relation, actually intertwines the coproduct of a Hopf
algebra with its opposite.
\vskip.5cm

\no {\bf Proposition 12.1.}  (a) There exists a unique homomorphism
$\Delta:~{\cal{A}}\to {\cal{A}}\otimes {\cal{A}}$, such that
$$
\eqalign{\Delta(H_a) &= H_a\otimes 1+1\otimes H_a~, \quad a\in M~, \cr
\Delta(e_\alpha) &= 1\otimes e_\alpha + e_\alpha \otimes
e^{\varphi(\alpha,\cdot)}~, \cr
\Delta(e_{-\alpha}) &= e^{-\varphi(\cdot,\alpha)} \otimes
e_{-\alpha}+e_{-\alpha}\otimes 1~, \quad \alpha \in N~. \cr}
\eqno(12.1)
$$
\no (b) If $I\subset {\cal{A}}$ is the ideal generated by the constants
in ${\cal{A}}^+$ and ${\cal{A}}^-$, and ${\cal{A}}^\prime= {\cal{A}}/I$,
then $\Delta$ induces a unique homomorphism
${\cal{A}}^\prime \to {\cal{A}}^\prime\otimes {\cal{A}}^\prime$ that will
also be denoted $\Delta$, so that (12.1) holds with $H_a$ and
$e_{\pm\alpha}$ being interpreted as generators of ${\cal{A}}/I$.
\s
\no (c) Let $\Delta^\prime$ be the opposite coproduct on ${\cal{A}}/I$,
and $R$ the standard R-matrix on ${\cal{A}}/I$ (satisfying Yang-Baxter),
then $\Delta R=R\Delta^\prime$.
\s
\no (d) The algebra ${\cal{A}}$ becomes a Hopf algebra when endowed with
the counit ${\cal E}$ and the antipode $S$.  The former is the unique
homomorphism ${\cal{A}}\to \Crm$ such that
$$
\eqalign{& {\cal E}(a)=1~, \quad {\cal E} (H_a)=0~, \quad a\in M~, \cr &
{\cal E}(e_{\pm\alpha})=0~, \quad \alpha\in N~. \cr}
\eqno(12.2)
$$
\no The antipode is the unique anti-automorphism $S:~{\cal{A}}\to
{\cal{A}}$ such that
$$
\eqalign{& S(1)=1~, \quad S(H_a)=-H_a~, \quad a\in M~, \cr &
S(e_\alpha)=-e_\alpha e^{-\varphi(\alpha,\cdot)}~, \quad
S(e_{-\alpha})=-e^{\varphi(\cdot,\alpha)} e_{-\alpha}~, \quad
\alpha\in N~. \cr} \eqno(12.3)
$$
\no (e) The counit ${\cal E}$ and the antipode $S$ of ${\cal{A}}$ induce
analogous structures on ${\cal{A}}^\prime={\cal{A}}/I$ such that (12.3)
holds on ${\cal{A}}^\prime$.
\vskip.5cm

\line  {\bf Proof. \hfil}
\s
 (a) The verification amounts to checking that $\Delta({\cal{A}})$ has
the relations of ${\cal{A}}$, in particular,
$$
\bigl[\Delta(e_\alpha),\Delta(e_{-\beta})\bigr] =
\delta^\beta_\alpha~\Delta\bigl([e_\alpha,e_{-\beta}]\bigr)~.
$$ (b) The ideal $I$ is generated by elements $x\in {\cal{A}}^+$ and
$y\in {\cal{A}}^-$ such that
$[e_{-\alpha},x]=0=[e_\alpha,y],~\alpha\in N$.  Since
$\Delta:~{\cal{A}}\to {\cal{A}}\otimes {\cal{A}}$ is a homomorphism,
$\Delta$ induces a homomorphism ${\cal{A}}/I\to ({\cal{A}}\otimes
{\cal{A}})/\Delta(I)$.  We must show that
$\Delta(I)\subset I\otimes {\cal{A}} + {\cal{A}}\otimes I$.  Since $I$ is
generated by elementary constants, it is enough to show that, for an
elementary constant $C$, $\Delta (C) \subset I\otimes {\cal{A}} +
{\cal{A}}\otimes I$.  Let $C\in {\cal{A}}^+$ be an elementary constant;
then $[e_{-\alpha},C]=0$ and thus $[\Delta(e_{-\alpha}),\Delta (C)]=0,~
\alpha\in N$.  If $C$ is of order $n$ in the generators, then (12.1)
shows that
$$
\Delta C=1\otimes C+ P^1\otimes P_{n-1}+P^2\otimes P_{n-2}+\ldots +
C\otimes P_0~,
$$
\no where $P^n$ and $P_k$ are homogeneous of order $k$ in the
$e_\alpha$'s.  Because $C$ is an elementary constant---Definition 4.1.---
there is no constant among the $P^k,P_k,~n=1,\ldots,n-1$; then
$[\Delta (e_{-\alpha}),\Delta (C)]=0$ implies that $\Delta C= 1\otimes
C+C\otimes P_0$ which indeed belongs to $ I\otimes {\cal{A}} +
{\cal{A}}\otimes I$; consequently $\Delta$ provides a map ${\cal{A}}/I\to
{\cal{A}}/I\otimes {\cal{A}}/I$.

(c) We use the abbreviation~--~ compare (10.6), Definition 2.2 and Eq.
(2.9),
$$ R=t_{(\alpha)}^{(\alpha^\prime)} R^i[e_{-\alpha}]\otimes
R_i[e_{\alpha^\prime}]~,
$$
$$
\eqalign{\Delta (e_\beta)R-R\Delta^\prime (e_\beta) &=
t^{(\alpha^\prime)}_{(\alpha)}
\bigl(R^i[e_{-\alpha}]\otimes e_\beta R_i[e_{\alpha^\prime}] \cr
&+e_\beta R^i[e_{-\alpha}]\otimes e^{\varphi(\beta,\cdot)}
R_i[e_{\alpha^\prime}]-R^i[e_{-\alpha}]e_\beta\otimes
R_i[e_{\alpha^\prime}] \cr
&-R^i[e_{-\alpha}]e^{\varphi(\beta,\cdot)}\otimes R_i[e_{\alpha^\prime}]
e_\beta\bigr)~. \cr}
$$
\no Terms 2 and 3 combine to $R^i[e_\beta,t^{(\alpha^\prime)}]\otimes
R_i[e_{\alpha^\prime}]$, and the recursion relations (7.2) implies that
the sum of all four terms is equal to zero.

(d) The existence and the uniqueness of the homomorphism ${\cal E}$ and
the anti-homomor-phism $S$ are obvious.  We have to show that
${\cal E}$ satisfies the axioms
$$ ({\cal E}\times id)\Delta=id=(id\times{\cal E})\Delta~,
$$
\no which is straightforward, and that
$$ m(id\times S)\Delta=\epsilon=m(S\times id)~.
$$
\no Here $m$ indicates multiplication, ${\cal{A}}\otimes {\cal{A}}\to
{\cal{A}}$.  For example,
$$
 m(id\times S)\Delta (e_\alpha)=S(e_\alpha)+e_\alpha
e^{-\varphi(\alpha,\cdot)}=0~.
$$ (e) Obvious, since ${\cal E}(I)=0$ and $S(I)=I$ by Proposition 7.2.
Theorem 12 is proved.

We turn to the case of the deformed R-matrix of Section 10, all
statements should be understood to hold to first order in the deformation
parameter $\epsilon$.  The maps $\Delta,{\cal E}$ and $S$ are    as
before and the deformed maps are
$
\Delta_\epsilon=\Delta + \epsilon\Delta_1~, \quad {\cal E}_\epsilon={\cal
E}+{\epsilon}{\cal E}_1~, \quad S_\epsilon=S+\epsilon S_1~.
$
\vskip.5cm
\no {\bf Proposition 12.2.}  (a) There is a unique homomorphism
$\Delta_\epsilon:~{\cal{A}}\to {\cal{A}}\otimes {\cal{A}}$ such that
$$
\Delta_1(x)=[\Delta(x),~Ke_{-\rho}\otimes Ke_\sigma]~, \quad x\in
{\cal{A}}~.
$$
\no (b) The projection of $\Delta_\epsilon$ to ${\cal{A}}^\prime\to
{\cal{A}}^\prime \otimes {\cal{A}}^\prime$ is well defined. (c) Let
$\Delta^\prime_\epsilon$ be the opposite coproduct on
${\cal{A}}^\prime={\cal{A}}/I$, and $R_\epsilon = R+\epsilon R_1$ the
R-matrix of Theorem 10, then $\Delta_\epsilon R_\epsilon= R_\epsilon
\Delta^\prime_\epsilon$ (to first order in $\epsilon$).
 (d) The deformed counit and antipode of ${\cal{A}}$ are given by
${\cal E}_1=0$ and
$$ S_1(x)=[Ke_{-\rho}e_\sigma,~S(x)]~, \quad x\in {\cal{A}}~.
$$
\no (e) The counit ${\cal E}_\epsilon$ and the antipode $S_\epsilon$
induce analogous structures on ${\cal{A}}/I$.
\vskip.5cm

\no {\bf Proof.}  (a) By the Jacobi identity.  (b) Obvious, for
$\Delta_1(C)=[\Delta(C),~Ke_{-\rho}\otimes Ke_\sigma]\in I\otimes
{\cal{A}}+{\cal{A}}\otimes I$.  (c) Completely straightforward. (d) We
have $({\cal E}\times id)\Delta_1(x)=0$, whence
${\cal E}_1=0$, while
$$ m(id\times S_1)\Delta(H_a) + m(id\times S) \Delta_1(H_a)=0
$$
\no since
$$  m(id\times S_1)\Delta(H_a)=S_1(H_a)= [H_a,~Ke_{-\rho}e_\sigma]~,
$$
$$
\eqalign{m(id\times S)\Delta_1(H_a) &= m(id\times S)(H_a(\sigma)-
H_a(\rho) Ke_{-\rho}\otimes Ke_\sigma \cr &=
\bigl(H_a(\sigma)-H_a(\rho)\bigr) Ke_{-\rho}\bigl(-e_\sigma
e^{-\varphi(\sigma,\cdot)}K^{-1}\bigr) \cr &=
-\bigl(H_a(\sigma)-H_a(\rho)\bigr) Ke_{-\rho}e_\sigma = -[H_a,Ke_{-\rho}
e_\sigma] \cr}
$$
\no and
$$ m(id\times S_1)\Delta(e_\alpha)+m(id\times S)\Delta_1(e_\alpha)=0
$$
\no since
$$
\eqalign{m(id\times S_1)\Delta(e_\alpha) &= S_1(e_\alpha)+e_\alpha
S_1(e^{\varphi(\alpha,\cdot)})\cr &= S_1(e_\alpha)-e_\alpha
\bigl[ e^{-\varphi(\alpha,\cdot)},~Ke_{-\rho}e_\sigma\bigr] \cr &=
[e_\alpha, Ke_{-\rho}e_\sigma]e^{-\varphi(\alpha,\cdot)}~. \cr m(id
\times S) \Delta_1(e_\alpha) &= -
[e_\alpha,Ke_{-\rho}e_\sigma]e^{-\varphi(\alpha,\cdot)}.\cr}
$$
\no These last two results require some work.
\s (e) This is clear, since ${\cal E}_1=0$ and $S_1(I)\in I$. The theorem
is proved.
\vskip1in

\parindent0pt {\bf References.}

\noindent

[BD]  A.A. Belavin and V.G. Drinfeld, Sov.Sci.Rev.Math.{\bf 4} (1984)
93-165.

[C]   C. Chevalley, ``Sur la classification des algebres de Lie simples
et de leur representations. \line {\hfil   C.R. {\bf 227} (1948)
1136-1138.}

[D1]  V.G. Drinfeld, in Proceedings of the International Congress of
Mathematicians, Berke- \line {\hfil ley  A.M. Gleason, ed. (American
Mathematical Society, Providence RI 1987.)}

[D2]  V.G. Drinfeld, in "Quantum Groups," P.P. Kulish, ed. Proceedings of
the Workshop \line { \hfil held in the Euler International  Mathematical
Institute, Leningrad, Fall 1990.}

[CG]~  E. Cremmer and J.-L. Gervais, Commun.Math.Phys. {\bf 134} (1990)
619.

[F]~~ ~  C. Fr\o nsdal, in Proceedings of the Karpac Winter School,
February 1994.

[FG1] C. Fr\o nsdal and A. Galindo, Contemporary Mathematics, {\bf 175}
(1994) 73-88.

[FG2] C. Fr\o nsdal and A. Galindo, Lett.Math.Phys. {\bf 34} (1995) 25-36.

[K] ~~  V.G. Kac, ``Infinite Dimensional Lie Algebras," Cambridge
University Press 1990.

[KR]~  A.N. Kirillov and N.Yu. Reshetikhin, Commun.Math.Phys.  {\bf 134}
(1991) 421-431.

[Ma]~  Yu.I. Manin, "Topics in noncommutative geometry," Princeton
University Press, \line {\hfil Princeton NJ 1991.}

[M] ~~ R.V. Moody, "A new class of Lie algebras." J. Algebra {\bf 10}
(1968) 211-230.

[LSo] J. Geom. Phys. {\bf 7} (1991) 1-14.

[R]~~ ~  N.Yu. Reshetikhin, Lett.Math.Phys. {\bf 20} (1990) 331-336.

[Ro]~~  M. Rosso, Commun.Math.Phys.{\bf 124} (1989) 307-319.

[Sc] ~ A. Schirrmacher, Z. Phys.C {\bf 50} (1991) 321.

[Su]~~  A. Sudbery, J.Phys.A.Math.Gen. {\bf 23} (1990) L697.

[V] ~~  V.S. Varadarajan

\bye